\documentclass[11pt]{article}
\usepackage[margin=1in]{geometry}
\usepackage{graphicx}
\usepackage{adjustbox} 
\usepackage[english]{babel}
\usepackage[sort, numbers]{natbib}

\usepackage{amsthm}
\usepackage{amssymb}
\usepackage{subcaption}
\usepackage[latin1]{inputenc}
\usepackage{enumerate}
\usepackage{amsmath}
\usepackage{hyperref}
\usepackage{color}
\usepackage{tikz}
\usepackage{pdflscape}
\usepackage{comment}
\usepackage{booktabs}
\usepackage[shortlabels]{enumitem}
\usepackage{xfrac}
\usepackage{float}
\usepackage[boxruled,vlined]{algorithm2e}
\usepackage[mathscr]{eucal}

\usepackage{siunitx}
\usepackage{booktabs}
\usepackage{multirow}
\usepackage{makecell}

\usepackage[
singlelinecheck=false 
]{caption}

\usepackage{pgfplots}
\usepackage{pgfplotstable}
\SetVlineSkip{0pt} 

\definecolor{darkblue}{rgb}{0.0,0.0,0.3}
\hypersetup{colorlinks,breaklinks,linkcolor=darkblue,urlcolor=darkblue,anchorcolor=darkblue,citecolor=darkblue}

\definecolor{MyGrey}{RGB}{100,120,100}
\definecolor{MyYellow}{RGB}{218,165,32}
\definecolor{MyGreen}{RGB}{0,102,51}
\definecolor{MyPurple}{RGB}{102,10,153}

\usetikzlibrary{shapes}
\usetikzlibrary{backgrounds,calc}
\usetikzlibrary{arrows}

\newtheorem{definition}{Definition}

\pgfplotsset{compat = 1.10}
\pgfplotsset{
        legend style={font=\scriptsize,row sep=-0.1cm,/tikz/every odd column/.append style={column sep=0.01cm}}
    }
\begin{document}

\title{A Robust Optimization Framework for Two-Echelon Vehicle and UAV Routing for Post-Disaster Humanitarian Logistics Operations}

\author{
	Tasnim Ibn Faiz\\ 
	 Department of Mechanical and Industrial Engineering \\
	Northeastern University \\ 
    Boston, MA 02115, USA \\
	\\
	Chrysafis Vogiatzis\\
 Industrial and Enterprise Systems Engineering \\
	University of Illinois at Urbana-Champaign \\ 
    Urbana, IL 61801, USA \\
    \\
	Md.\ Noor-E-Alam\footnote{Corresponding author. Email: mnalam@neu.edu.} \\
 Department of Mechanical and Industrial Engineering \\
	Northeastern University \\ 
    Boston, MA 02115, USA
}
\date{}

\maketitle

\begin{abstract}

Providing first aid and other supplies (e.g., epi-pens, medical supplies, dry food, water) during and after a disaster is always challenging. The complexity of these operations increases when the transportation, power, and communications networks fail, leaving people stranded and unable to communicate their locations and needs. The advent of emerging technologies like uncrewed autonomous vehicles can help humanitarian logistics providers reach otherwise stranded populations after transportation network failures. However, due to the failures in telecommunication infrastructure, demand for emergency aid can become uncertain. To address the challenges of delivering emergency aid to trapped populations with failing infrastructure networks, we propose a novel robust computational framework for a two-echelon vehicle routing problem that uses uncrewed autonomous vehicles, or drones, for the deliveries. We formulate the problem as a two-stage robust optimization model to handle demand uncertainty. Then, we propose a column-and-constraint generation approach for worst-case demand scenario generation for a given set of truck and drone routes.
Moreover, we develop a decomposition scheme inspired by the column generation approach to heuristically generate drone routes for a set of demand scenarios. Finally, we combine the heuristic decomposition scheme within the column-and-constraint generation approach to determine robust routes for both trucks and drones, the time that affected communities are served, and the quantities of aid materials delivered. To validate our proposed computational framework, we use a simulated dataset that aims to recreate emergency aid requests in different areas of Puerto Rico after Hurricane Maria in 2017.

\noindent 

\end{abstract}

\section{Introduction}

In this work, we investigate a paradigm for humanitarian relief operations in the absence of stable transportation, communication, and power infrastructures. Providing relief to affected populations in the form of first aid items, blankets, water, and dry food, is difficult even in the presence of a fully functional built environment. When the infrastructure itself is prone to failures, these operations become even more challenging as people may become stranded without access to telecommunications, transportation, or power. The advent of uncrewed autonomous vehicles (UAVs) provides us with a reliable means of reaching populations otherwise unreachable due to the failures in the built environment. At the same time, we can do so with fewer risks for the first responders. While UAVs are naturally suited to improving safety under these conditions, they also present us with new challenges, namely the inherent computational difficulty of such routing operations and the uncertainty of the location and number of demand requests to which they need to respond. To that end, we propose a new robust optimization framework to incorporate demand information uncertainty within a two-echelon vehicle routing problem to aid trapped populations using ground vehicles and deliver it using UAVs.

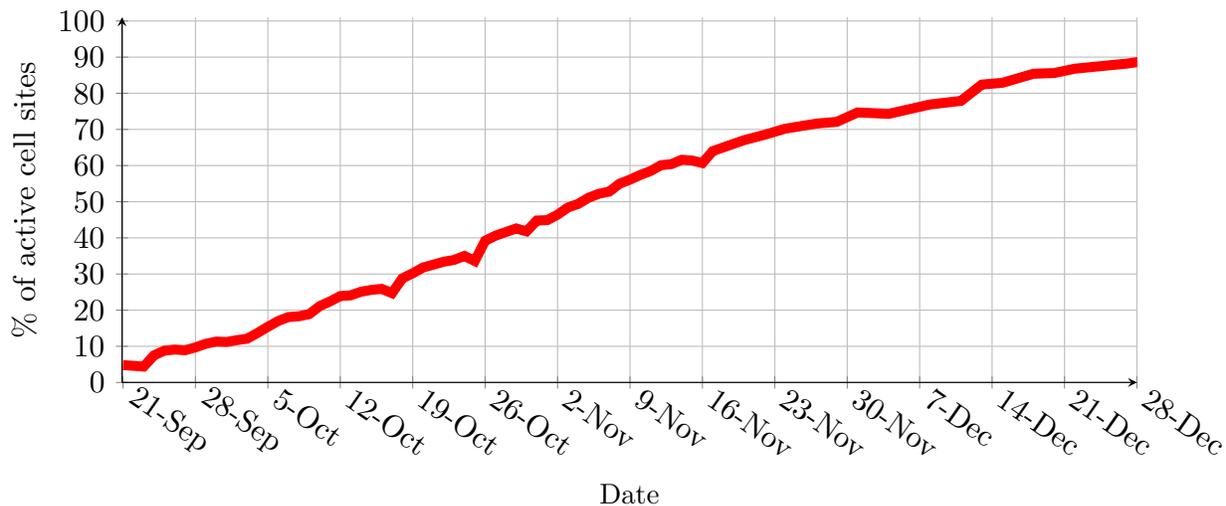
\begin{figure}[htbp]
\caption{Percentage of cell sites that were active and working to provide access to cell phone reception since landfall (September 21, 2017) and until January 4, 2018. More than 90 days after Hurricane Maria's landfall, the number was still shy of 100\%. More importantly, three weeks after landfall, that percentage was only 23.7\%. All data was taken by the Federal Communications Commission (FCC) \cite{FCC2018}. \label{fcc1}}
\centering
\resizebox{\textwidth}{!}{\begin{tikzpicture}
\begin{axis}[%
axis lines=left,
width=12.5cm,
height=4.5cm,
scale only axis,
xmin=-0.1,
xmax=98,
xtick={0, 7, 14, 21, 28, 35, 42, 49, 56, 63, 70, 77, 84, 91, 98},
xticklabels={21-Sep, 28-Sep, 5-Oct, 12-Oct, 19-Oct, 26-Oct, 2-Nov, 9-Nov, 16-Nov, 23-Nov, 30-Nov, 7-Dec, 14-Dec, 21-Dec, 28-Dec, 4-Jan},
ylabel={\% of active cell sites},
xlabel={\small Date},
x tick label style={rotate=325,anchor=west},
xmajorgrids=true,
ymajorgrids=true,
ymin=0,
ymax=101,
ytick={0, 10, 20, 30, 40, 50, 60, 70, 80, 90, 100}]

        \addplot[color=red, line width=3.5pt] coordinates {
(0,4.8)
(1,4.59999999999999)
(2,4.40000000000001)
(3,7.5)
(4,8.8)
(5,9.09999999999999)
(6,8.90000000000001)
(7,9.7)
(8,10.7)
(9,11.3)
(10,11.2)
(11,11.7)
(12,12.1)
(13,13.7)
(14,15.4)
(15,17)
(16,18.1)
(17,18.3)
(18,18.9)
(19,21.1)
(20,22.4)
(21,23.9)
(22,24.1)
(23,25.1)
(24,25.6)
(25,25.9)
(26,24.7)
(27,28.8)
(28,30.2)
(29,31.8)
(30,32.6)
(31,33.4)
(32,33.9)
(33,35)
(34,33.6)
(35,39.2)
(36,40.6)
(37,41.6)
(38,42.6)
(39,41.8)
(40,44.8)
(41,44.9)
(42,46.4)
(43,48.4)
(44,49.4)
(45,51.1)
(46,52.2)
(47,52.8)
(48,55)
(49,56.1)
(50,57.4)
(51,58.5)
(52,60.1)
(53,60.4)
(54,61.6)
(55,61.4)
(56,60.7)
(57,64)
(60,67)
(62,68.5)
(64,70.2)
(67,71.6)
(69,72.1)
(71,74.7)
(74,74.3)
(76,75.6)
(78,76.9)
(81,77.9)
(83,82.4)
(85,82.9)
(88,85.4)
(90,85.6)
(92,86.8)
(97,88.2)
(99,89)
(104,90.5)
(106,90.8)
(109,90.4)
};

\end{axis}
\end{tikzpicture}}%
\end{figure}

\begin{figure}[htbp]
\caption{Number of counties (out of a total of 78 on the island of Puerto Rico) with no active cell site since landfall and until November 9, 2017. All data was taken by the Federal Communications Commission (FCC) \citep{FCC2018}. \label{fcc2}}
\centering
 \resizebox{0.55\textwidth}{!}{\begin{tikzpicture}
    \begin{axis}[
    	scaled ticks=false, tick label style={/pgf/number format/fixed},
    	xmin=0, xmax=50,
        ymin=0, ymax=50,
        samples=400,
        axis y line=left,
        axis x line=bottom,
        xmajorgrids=true,
        ymajorgrids=true,
        xlabel=Days after landfall,
        ylabel=\# Puerto Rico counties 
    ]												
        \addplot[color=blue, line width=2pt] coordinates {
(0,48)
(1,47)
(2,48)
(3,37)
(4,34)
(5,29)
(6,31)
(7,29)
(8,29)
(9,28)
(10,27)
(11,27)
(12,33)
(13,27)
(14,24)
(15,22)
(16,22)
(17,22)
(18,23)
(19,18)
(20,15)
(21,14)
(22,12)
(23,12)
(24,11)
(25,11)
(26,12)
(27,8)
(28,8)
(29,5)
(30,4)
(31,4)
(32,4)
(33,4)
(34,4)
(35,4)
(36,4)
(37,1)
(38,0)
(39,0)
(40,0)
(41,0)
(42,1)
(43,0)
(44,0)
(45,0)
(46,0)
(47,0)
(48,0)
(49,0)
(50,0)
	};
\end{axis}
\end{tikzpicture}}
\end{figure}
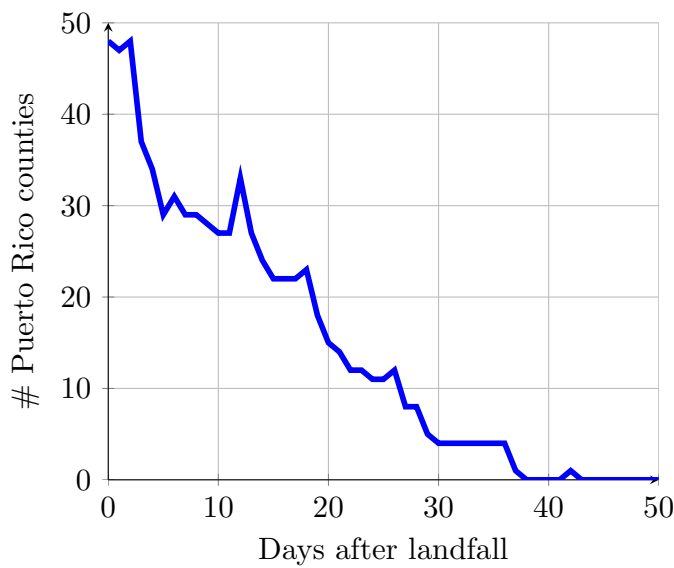

We are primarily motivated by the situation in Puerto Rico after Hurricane Maria in 2017. The hurricane made landfall on September 21, and eventually left behind an island without much of its built environment: power, and consequently communications, were unavailable even as much as three weeks after the disaster hit (see Figures \ref{fcc1} and \ref{fcc2}, as well as newspaper articles, e.g., \cite{WaPo2017}). All 78 counties in Puerto Rico had more than 75\% of their cell sites out of service until September 29, 2017. Even then, most counties had no cell phone coverage, with 100\% of their cell sites being deactivated during and after landfall, as shown in Figure \ref{fcc2}. It took 44 days for all counties to have at least some cellular coverage. On a similar note, 46 days passed until every county had at least half its cell sites reinstated. In addition, people were also stranded due to damaged roads, making it impossible to know their locations and conditions (no communication) or drive to them (no transportation) for days \cite{LATimes2017}. Our proposed framework would help provide assistance and valuable items to these trapped populations in a robust manner.

The paper is outlined as follows. First, we provide a brief literature review in Section \ref{sec:literature}. In Section \ref{sec:formulation}, we proceed to define our problem formally and provide our formulation, followed by our novel solution approach. In that same section, we describe all aspects of our robust framework, including the proposed column generation approach and the scenario generation subproblem. Section \ref{sec:results} provides a case study on a simulated dataset that aims to capture the demands in Puerto Rico after Hurricane Maria. We conclude our work in Section \ref{sec:conclusions} with a summary of our findings and plans for future work in this area.

\section{Literature review and contributions}
\label{sec:literature}

In this section, we provide a brief literature review, before we present our main contributions. 

\subsection{Review of relevant research}

We begin with the formal definition of the two-echelon vehicle routing problem (2EVRP) \cite{crainic2009models,perboli2011two,cuda2015survey}, an interesting extension of the traditional vehicle routing problem. 

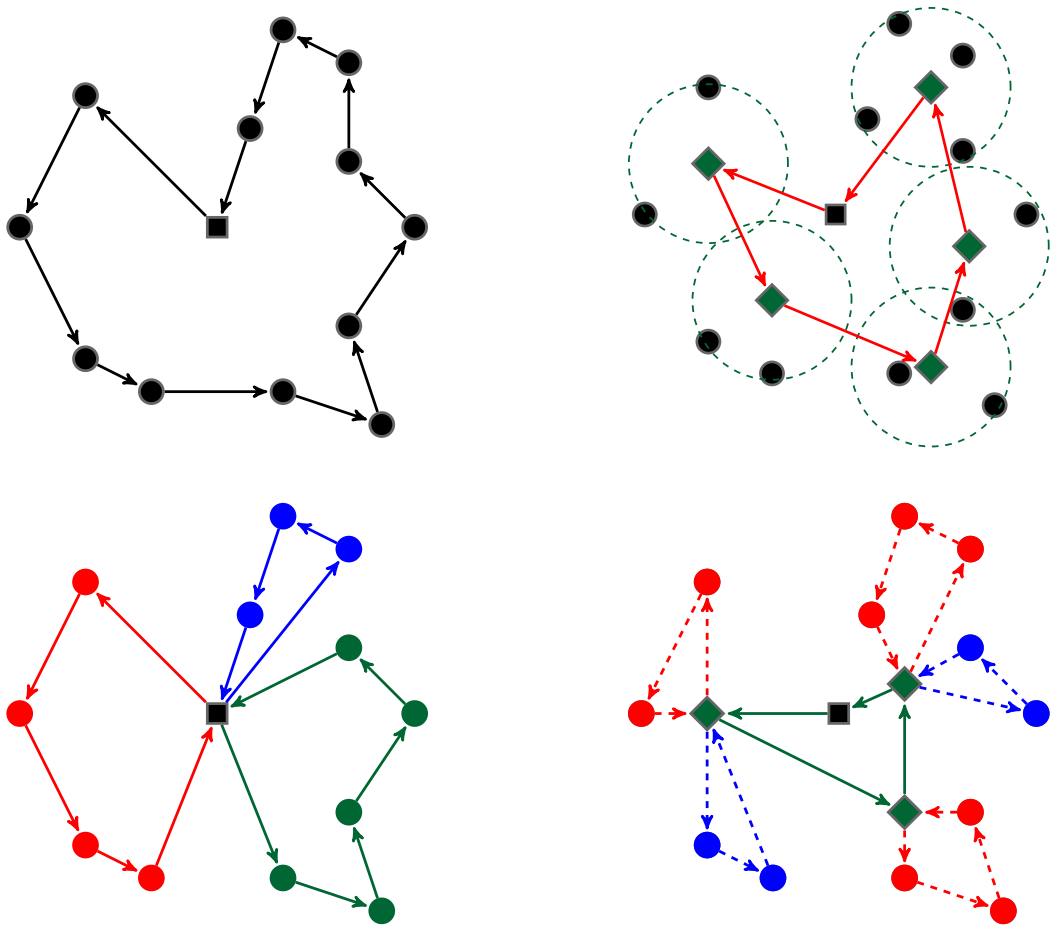
\begin{figure}[t]
\centering
\caption{Four related problems in the literature. The Traveling Salesperson Problem (TSP, top left), the Covering Traveling Salesperson Problem (Covering TSP, top right), the traditional Vehicle Routing Problem (VRP, bottom left), and the Two-Echelon Vehicle Routing Problem (2EVRP, bottom right). \label{fig:graphical_literature_review}}
	\begin{minipage}{.485\textwidth}
    \centering
    	\resizebox{.7\textwidth}{!}{\begin{tikzpicture}[
roundnode/.style={circle, draw=black!60, fill=black, very thick, minimum size=1.5mm},
squarednode/.style={rectangle, draw=black!60, fill=black, very thick, minimum size=3mm},
->,>=stealth',shorten >=1pt,auto,node distance=1.5cm,semithick]

\node[squarednode]  at (1,0) (depot)    {};
\node[roundnode]  at (-1,2)  (1)            {};
\node[roundnode]  at (-2,0)  (2)            {};
\node[roundnode]  at (-1,-2)  (3)            {};
\node[roundnode]  at (0,-2.5)  (4)            {};
\node[roundnode]  at (2,-2.5)  (5)            {};
\node[roundnode]  at (3.5,-3)  (6)            {};
\node[roundnode]  at (3,-1.5)  (7)            {};
\node[roundnode]  at (4,0)  (8)            {};
\node[roundnode]  at (3,1)  (9)            {};
\node[roundnode]  at (3,2.5)  (10)            {};
\node[roundnode]  at (2,3)  (11)            {};
\node[roundnode]  at (1.5,1.5)  (12)            {};

\path (depot) edge[very thick]                   (1);
\path (1) edge[very thick]                     (2);
\path (2) edge[very thick]      			(3);
\path[bend right=0] (3) edge[very thick]      (4);
\path (4) edge[very thick]      			(5);
\path (5) edge[very thick]      			(6);
\path (6) edge[very thick]      			(7);
\path (7) edge[very thick]      			(8);
\path (8) edge[very thick]      			(9);
\path (9) edge[very thick]      			(10);
\path (10) edge[very thick]      			(11);
\path (11) edge[very thick]      			(12);
\path (12) edge[very thick]      			(depot);

\end{tikzpicture}}
\end{minipage}~
\begin{minipage}{.485\textwidth}
\centering
    	\resizebox{.7\textwidth}{!}{\begin{tikzpicture}[
roundnode/.style={circle, draw=black!60, fill=black, very thick, minimum size=1.5mm},
squarednode/.style={rectangle, draw=black!60, fill=black, very thick, minimum size=3mm},
trianglenode/.style={diamond, draw=black!60, fill=MyGreen, very thick, minimum size=2.25mm},
->,>=stealth',shorten >=1pt,auto,node distance=1.5cm,semithick]

\node[squarednode]  at (1,0) (depot)    {};
\node[roundnode]  at (-1,2)  (1)            {};
\node[roundnode]  at (-2,0)  (2)            {};
\node[roundnode]  at (-1,-2)  (3)            {};
\node[roundnode]  at (0,-2.5)  (4)            {};
\node[roundnode]  at (2,-2.5)  (5)            {};
\node[roundnode]  at (3.5,-3)  (6)            {};
\node[roundnode]  at (3,-1.5)  (7)            {};
\node[roundnode]  at (4,0)  (8)            {};
\node[roundnode]  at (3,1)  (9)            {};
\node[roundnode]  at (3,2.5)  (10)            {};
\node[roundnode]  at (2,3)  (11)            {};
\node[roundnode]  at (1.5,1.5)  (12)            {};

\node[trianglenode] at (2.5,2) (center5) {};
\draw[MyGreen,thick,dashed] (center5) circle (1.25);

\node[trianglenode] at (3.1,-0.5) (center4) {};
\draw[MyGreen,thick,dashed] (center4) circle (1.25);

\node[trianglenode] at (2.5,-2.4) (center3) {};
\draw[MyGreen,thick,dashed] (center3) circle (1.25);

\node[trianglenode] at (0,-1.35) (center2) {};
\draw[MyGreen,thick,dashed] (center2) circle (1.25);

\node[trianglenode] at (-1,0.8) (center1) {};
\draw[MyGreen,thick,dashed] (center1) circle (1.25);


\path[red] (depot) edge[very thick]                   (center1);
\path[red] (center1) edge[very thick]                     (center2);
\path[red] (center2) edge[very thick]      			(center3);
\path[red,bend right=0] (center3) edge[very thick]      (center4);
\path[red] (center4) edge[very thick]      			(center5);
\path[red] (center5) edge[very thick]      			(depot);

	\end{tikzpicture}}
\end{minipage} 

\vspace{20pt}

\begin{minipage}{.485\textwidth}
\centering
\resizebox{.7\textwidth}{!}{\begin{tikzpicture}[
roundnode/.style={circle, draw=black!60, fill=black, very thick, minimum size=1.5mm},
squarednode/.style={rectangle, draw=black!60, fill=black, very thick, minimum size=3mm},
->,>=stealth',shorten >=1pt,auto,node distance=1.5cm,semithick]

\node[squarednode]  at (1,0) (depot)    {};
\node[roundnode, red]  at (-1,2)  (1)            {};
\node[roundnode, red]  at (-2,0)  (2)            {};
\node[roundnode, red]  at (-1,-2)  (3)            {};
\node[roundnode, red]  at (0,-2.5)  (4)            {};
\node[roundnode, MyGreen]  at (2,-2.5)  (5)            {};
\node[roundnode, MyGreen]  at (3.5,-3)  (6)            {};
\node[roundnode, MyGreen]  at (3,-1.5)  (7)            {};
\node[roundnode, MyGreen]  at (4,0)  (8)            {};
\node[roundnode, MyGreen]  at (3,1)  (9)            {};
\node[roundnode, blue]  at (3,2.5)  (10)            {};
\node[roundnode, blue]  at (2,3)  (11)            {};
\node[roundnode, blue]  at (1.5,1.5)  (12)            {};

\path[red] (depot) edge[very thick]                   (1);
\path[red] (1) edge[very thick]                     (2);
\path[red] (2) edge[very thick]      			(3);
\path[red,bend right=0,very thick] (3) edge[very thick]      (4);
\path[red,bend right=0] (4) edge[very thick]      (depot);

\path[MyGreen] (depot) edge[very thick]      		(5);
\path[MyGreen] (5) edge[very thick]      			(6);
\path[MyGreen] (6) edge[very thick]      			(7);
\path[MyGreen] (7) edge[very thick]      			(8);
\path[MyGreen] (8) edge[very thick]      			(9);
\path[MyGreen] (9) edge[very thick]      			(depot);

\path[blue] (depot) edge[very thick]      		(10);
\path[blue] (10) edge[very thick]      			(11);
\path[blue] (11) edge[very thick]      			(12);
\path[blue] (12) edge[very thick]      			(depot);

\end{tikzpicture}}
\end{minipage}~
\begin{minipage}{.485\textwidth}
\centering
\resizebox{.7\textwidth}{!}{\begin{tikzpicture}[
roundnode/.style={circle, draw=black!60, fill=black, very thick, minimum size=1.5mm},
squarednode/.style={rectangle, draw=black!60, fill=black, very thick, minimum size=3mm},
trianglenode/.style={diamond, draw=black!60, fill=MyGreen, very thick, minimum size=2.25mm},
->,>=stealth',shorten >=1pt,auto,node distance=1.5cm,semithick]

\node[squarednode]  at (1,0) (depot)    {};

\node[roundnode]  at (-1,2)  (1)            {};
\node[roundnode]  at (-2,0)  (2)            {};
\node[roundnode]  at (-1,-2)  (3)            {};
\node[roundnode]  at (0,-2.5)  (4)            {};
\node[roundnode]  at (2,-2.5)  (5)            {};
\node[roundnode]  at (3.5,-3)  (6)            {};
\node[roundnode]  at (3,-1.5)  (7)            {};
\node[roundnode]  at (4,0)  (8)            {};
\node[roundnode]  at (3,1)  (9)            {};
\node[roundnode]  at (3,2.5)  (10)            {};
\node[roundnode]  at (2,3)  (11)            {};
\node[roundnode]  at (1.5,1.5)  (12)            {};

\node[roundnode, red]  at (-1,2)  (1)            {};
\node[roundnode, red]  at (-2,0)  (2)            {};
\node[roundnode, blue]  at (-1,-2)  (3)            {};
\node[roundnode, blue]  at (0,-2.5)  (4)            {};
\node[roundnode, red]  at (2,-2.5)  (5)            {};
\node[roundnode, red]  at (3.5,-3)  (6)            {};
\node[roundnode, red]  at (3,-1.5)  (7)            {};
\node[roundnode, blue]  at (4,0)  (8)            {};
\node[roundnode, blue]  at (3,1)  (9)            {};
\node[roundnode, red]  at (3,2.5)  (10)            {};
\node[roundnode, red]  at (2,3)  (11)            {};
\node[roundnode, red]  at (1.5,1.5)  (12)            {};

\node[MyGreen,trianglenode] at (-1,0) (center1) {};
\node[MyGreen,trianglenode] at (2,-1.5) (center2) {};
\node[MyGreen,trianglenode] at (2,0.45) (center3) {};


\path[red] (center1) edge[very thick,dashed]                     (1);
\path[red] (1) edge[very thick,dashed]      			(2);
\path[red] (2) edge[very thick,dashed]      			(center1);

\path[blue] (center1) edge[very thick,dashed]                     (3);
\path[blue] (3) edge[very thick,dashed]      			(4);
\path[blue] (4) edge[very thick,dashed]      			(center1);

\path[red] (center2) edge[very thick,dashed]      		(5);
\path[red] (5) edge[very thick,dashed]      			(6);
\path[red] (6) edge[very thick,dashed]      			(7);
\path[red] (7) edge[very thick,dashed]      			(center2);

\path[blue] (center3) edge[very thick,dashed]      		(8);
\path[blue] (8) edge[very thick,dashed]      			(9);
\path[blue] (9) edge[very thick,dashed]      			(center3);

\path[red] (center3) edge[very thick,dashed]      		(10);
\path[red] (10) edge[very thick,dashed]      			(11);
\path[red] (11) edge[very thick,dashed]      			(12);
\path[red] (12) edge[very thick,dashed]      			(center3);

\path[MyGreen] (depot) edge[very thick]      		(center1);
\path[MyGreen] (center1) edge[very thick]      		(center2);
\path[MyGreen] (center2) edge[very thick]      		(center3);
\path[MyGreen] (center3) edge[very thick]      		(depot);

\end{tikzpicture}}
\end{minipage}
\end{figure}

\begin{definition}[2EVRP]
		Given a \underline{graph $G(V,E)$}, where the nodeset $V$ consists of $0$ (depot), a set of candidate locations $S$, and a set of customers $C$ (i.e., $V=\{0\}\cup S\cup C$), \underline{real values $\tau_{ij}$} for the costs of every edge $(i,j)\in E$, and \underline{real values $f_i$} for the costs of using a node $i$, the 2EVRP aims to identify: (a) the set of facilities $\mathcal{S}\subseteq S$ to use; (b) the circuit from $0$ to every node $s\in\mathcal{S}$; (c) the assignment of customers $C$ to a node $s\in\mathcal{S}$; and (d) the circuits from $s\in\mathcal{S}$ to every customer assigned to $s$ with minimum total cost. 
\end{definition} We note here that the problem can be easily shown to be $\mathcal{NP}$-hard, by transforming it to the traditional vehicle routing problem by letting only one satellite station. The 2EVRP is an extension of the traditional vehicle routing problem, where a set of customer locations needs to be visited by one of the tours that start and end at a facility. It is also very closely related to the covering traveling salesperson problem \cite{current1989covering,gendreau1997covering,golden2012generalized} and its ``close-enough'' variants \cite{gulczynski2006close,behdani2014integer}, where we identify a circuit that visits a subset of the nodes such that every other node is in close proximity to at least one visited node. We present the relationship and similarities between some of these two problems in visual form in Figure \ref{fig:graphical_literature_review}.

Exact solution methods for the 2EVRP are primarily based on branch-and-cut algorithms, such as the one proposed in \cite{gonzalez2008two}, which later served as a starting point for numerous follow-up studies \cite{perboli2010new,jepsen2013branch,perboli2011two}. Later, Santos et al.\ used column generation in the setting of 2EVRPs to propose both a branch-and-price approach \cite{santos2013branch} and a branch-and-cut-and-price \cite{santos2015branch} with additional valid inequalities. Strong valid inequalities have also been recently proposed in \cite{liu2018branch,perboli2018new}. As expected, heuristic approaches have been put to the use, as is the case with the large neighborhood search procedures proposed in \cite{hemmelmayr2012adaptive,breunig2016large}.

As can be seen by the setup of the 2EVRP, it offers a suitable framework to model new logistics technologies, such as the use of uncrewed autonomous vehicles for last-mile deliveries. Due to the advent of UAVs, research in vehicle routing problems with the use of drones is booming, as is obvious by the ever increasing number of publications and scholarly activities involving UAVs in logistics. One of the first attempts at mathematically modeling a coordination problem between a truck and a drone is made in \cite{murray2015flying}, where a traditional vehicle and its ``flying sidekick'' (a UAV/drone) satisfy static and deterministic demand in a geographic location; an extension to multiple ``sidekicks'' has appeared recently in \cite{dell2021modeling}. Later, the traveling salesperson problem with drones and the vehicle routing problem with drones were introduced in \cite{agatz2018optimization} and \cite{wang2017vehicle,poikonen2017vehicle}, respectively. Branch-and-cut, which has gained popularity in solving 2EVRPs, has also been shown to work well in practice with the TSP with drones \cite{schermer2020b}. Dynamic programming has also been successful in dealing with the TSP with drones problem, leading to faster solutions in larger instances compared to using commercial solvers \cite{bouman2018dynamic}.

Interestingly, due to the nature of the logistics problems that are usually of interest, heuristics involving route-first, cluster/assign-later approaches are very prominent. In \cite{agatz2018optimization}, the authors propose route-first cluster-second heuristics and derive approximation guarantees for their performance, while in \cite{POIKONEN2020104802}, a TSP solution is transformed eventually into a solution involving drone deliveries in a smart, flexible manner. The idea of transforming a TSP solution to a tour that involves drones is prominent in other recent works as well, either with the use of $k$-means clustering first \cite{CHANG2018307}, adaptive insertion \cite{KITJACHAROENCHAI201914}, or local searches and Greedy Randomized Adaptive Search Procedures (GRASP) \cite{HA2018597}. There also exist other heuristic procedures to help solve problems where a customer can be visited by either a truck or a drone, as in \cite{di2021trucks}.

Moreover, math-based heuristics that are motivated by well-known combinatorial optimization techniques, such as branch-and-bound or column generation, are also of interest to our community. For example, heuristics are employed in an attempt to prune the branch-and-bound tree and speed up computations in \cite{poikonen2019branch}. A different heuristic approach, inspired by column generation, is proposed and used in \cite{our_paper} to solve a problem with drones in a humanitarian relief setting, where part of the drones are tasked with providing cell phone coverage to an affected area in addition to the drones tasked with making deliveries. This application of using drones to deliver medical aid is also explored in \cite{GHELICHI2021105443}: the main difference to our work is that here, we offer a robust optimization approach to tackle the uncertainty of demand locations and demand quantities. Finally, a recent work dispels with the discretization of locations from where a drone can be deployed and/or received, and allows for drone movements in a continuous manner, that is, we are allowed to deploy and receive a drone on a moving truck while it is traversing an edge \cite{masone2022multivisit}.
 
Since our decision problem is stochastic in nature, we continue this literature review with the stochastic VRP. Uncertainty can arise in multiple aspects of the problem including demands, deadlines/time windows, costs/times, etc. Traditional vehicle routing problems with and without capacity constraints have addressed uncertainty in a number of ways. An early work in the area by Gendreau et al.\ proposed a two-stage stochastic integer programming framework, where routes are designed in the first stage and during the second stage the locations with no demand are skipped after demand realization \cite{gendreau1995exact}. In \cite{sungur2010model}, the authors adopt a robust optimization approach with scenario-based decomposition to solve a vehicle routing problem with stochastic service times. Probabilistic versions of the traveling salesperson problem have also been proposed over the years, as in the probabilistic traveling salesperson problem with deadlines \cite{campbell2008probabilistic} where customers are stochastic and, additionally, need to be visited before a certain time. An extension, where deadlines are replaced by time windows, is offered in \cite{voccia2013probabilistic}. A dynamic vehicle routing problem with time windows is proposed and then solved by column generation in \cite{chen2006dynamic}. Hvattum et al.\ consider the case when both customer locations and demands are only partially known in advance: that is, customers can call in orders with their demands and time windows. The problem is modeled as a dynamic and stochastic programming problem, and a multistage heuristic is proposed \cite{hvattum2006solving}. A similar dynamic problem with stochastic customers and partially available information a priori is solved in \cite{bent2004scenario}, with routes being generated continuously for the available scenarios, and a consensus methodology is adopted to select the best routes.

Several works adopting two-stage stochastic and robust optimization approaches for handling uncertainties in VRPs can be found in the literature. In contrast, literature on exact approaches for solving 2EVRP with stochastic parameters is very limited. Few available works have addressed such problems with  simulation-based \cite{liu2017simulation} and metaheuristic approaches (e.g., genetic algorithms \cite{wang2017genetic}). With our contributions in this work, we aim to fill this gap in the literature.

We finish this literature review with a brief discussion of recent works in vehicle routing problems in humanitarian aid and relief operations. Interested readers are referred to the excellent surveys in \cite{kovacs2007humanitarian,caunhye2012optimization,leiras2014literature,banomyong2019systematic}. Vehicle routing problems are, of course, a prominent aspect of humanitarian logistics especially in the aftermath of natural and human-made disasters (see, e.g., \cite{ozdamar2004emergency,ozdamar2012hierarchical,faiz2019column}). However, the underlying assumption of these studies is that post disaster vehicle routing and aid delivery happen in a fully functional transportation network, which unfortunately is less than realistic in many disasters. In the absence of a stable transportation network, emerging technologies (such as the use of UAVs) and innovative approaches (such as stochastic and robust frameworks) are necessary, and this is indeed the premise of our study. 

\subsection{Our contributions}

Our work focuses on addressing two critical challenges. First, from the side of humanitarian logistics, we assume that the infrastructure is unstable or unusable for ground vehicles after a disaster. UAVs are included in the humanitarian logistics planning operations to deliver medical items to trapped communities that are unreachable via ground transportation network. Secondly, since the telecommunication infrastructure has failed, we cannot assume that we have a priori knowledge of demand locations and quantities. The stochasticity in the parameters of modern routing problems involving UAVs has already been noted as an important future research direction \cite{poikonen2021future}. Therefore, we propose a novel robust optimization framework, that we refer to as R-2EVRP-UAV (\textbf{R}obust \textbf{2}-\textbf{E}chelon \textbf{V}ehicle \textbf{R}outing \textbf{P}roblem with \textbf{UAV}s) to tackle the unavailable information on the location and the quantities of demands. The framework is described in more detail in Section \ref{subsec:framework}. Specifically, our contributions can be summarized as follows:

\begin{itemize}
    \item We investigate the 2EVRP with stochastic demands in a post-disaster setting. In this setting, trucks can use part of the transportation infrastructure and can carry drones that are deployed towards communities to deliver aid. Our motivating assumption is that aid demands at the communities are not known in advance (uncertain demands) prior to route determination, hence the need for a robust optimization framework.
    \item Then, we present a decomposition framework and a heuristic approach, based on column generation techniques for generating UAV routes. Our decomposition approach leverages the problem structure and constructs small size  pricing subproblems, which enables us to solve large-scale problems.
    \item Finally, we develop a novel decision framework (R-2EVRP-UAV). For handling demand uncertainty, we design a column-and-constraint generation algorithm to generate demand scenarios. We then solve the decision model with an innovative computational framework that includes an outer level column-and-constraint generation and an inner level column generation method. Our model formulation and solution approach for robust 2EVRP with drones fills the gap in the state-of-the-art robust approach for handling uncertainty in 2EVRPs.
\end{itemize}

\section{Mathematical formulation}
\label{sec:formulation}

\subsection{Notation}

We begin with the necessary notation for our problem. Let $G(V,E)$ be a simple graph, representing the transportation network and the UAV aerial routes available. We specifically have that $V=\left\{0\right\}\cup S\cup C$, where $0$ is the depot, $S$ is an already identified set of candidate satellite locations, and $C$ is the set of affected communities. Note that there exist certain communities in $C$ that are reachable from a UAV departing from satellite station $s\in S$ (i.e., these communities are within the range of satellite station $s$). We note these communities as $C_s$ and, thus we have that $C=\bigcup_{s\in S} C_s$. 

Moreover, we let $E=E^t\cup E^d$, where $E^t$ is the set of edges that a truck can use; on the other hand, UAVs can use any edge in the network (hence, $E^d$ are the edges that can only be used by a UAV). We can also define subsets $E^d_s$, consisting of all drone edges that we can use when the truck is located at satellite station $s$, and hence we have $E^d=\bigcup_{s\in S}E^d_s$.

We make the assumption that all of the trucks and drones employed are identical. We also assume that weather conditions are deterministic and vehicle speeds (truck and drone) are constant and known. Hence, for all edges in the network, we associate constant traversal times denoted by $t$ for trucks and $d$ for UAVs. Specifically, $\tau_{ij}^t$ defined for all $\left(i,j\right)\in E^t$ are the traversal times for trucks and $\tau_{ij}^d$ defined for all $\left(i,j\right)\in E^d$ are the traversal times for drones. 

We use $D_s$ to represent the drones available at a truck that has stopped at satellite station $s$. Similarly, we define a set of available routes for a drone originating from satellite station $s$. As those are exponentially many, we only use a restricted set of routes, denoted as $V_s, \forall s\in S$. Note that if community $c\in C_s$ is visited during a route $p\in V_s$, then we write that $b_c^p=1$; $0$, otherwise. We also let $\Omega$ represent all scenarios, with $\Omega^\prime$ being the set of scenarios included in the main problem. $Q_c^\omega$ represents the quantity of items that community $c$ would request under scenario $\omega$.

For a given UAV $\ell\in D_s$ and route $p\in V_s$ originating from satellite station $s\in S$, we define $\Pi_{sp}^d$ as the total route duration, $\overline{t}_c^{s\ell}$ as the visiting time of community $c$, and $d_{sc}^{\omega p}$ as the delivery quantity to community $c$ for every scenario $\omega$. These quantities are all known, provided we have already generated the route (as will be shown in our pricing subproblem). 

Continuing with the parameters of our model, we have $m^t$ and $m^d$ be the number of trucks available and the number of drones available per truck, respectively. $F_c^T$ is the penalty cost per unit of time for delaying in aid delivery to community $c$; $F_c^R$ is the penalty cost for failing to reach a community altogether; $F_c^D$ is the shortage cost per unit of unmet demand. Finally, each drone has a flying range of $W^d$. We measure this range in ``units of time'': it can then be viewed as a budget for how much time the UAV can spend before it needs to return to the truck at the satellite station. A similar modeling assumption has been followed in previous work (see, e.g., \cite{doi:10.1287/ijoc.2018.0879,our_paper}). Drones are also restricted by a maximum carry load, denoted by $L_{max}$.

We are now ready to proceed with the description of our framework in subsection \ref{subsec:framework}, before moving to our main problem (in subsection \ref{subsec:main_problem}), our pricing subproblem (in subsection \ref{subsec:pricing}), and our scenario generation subproblem (in subsection \ref{subsec:scenario_generation}).

\subsection{Framework}
\label{subsec:framework}

\begin{figure}
    \centering
    \caption{A big-picture representation of the problem we are solving. One truck, starting from the depot, is allowed to use the transportation network (which connects the numbered nodes, in blue) in order to approach the affected communities (in red). Once the truck visits a node, it is allowed to deploy two drones (red and green) which visit a subset of the communities each, before returning to the truck. A yellow shaded circle with a number node (satellite station) as its center and dash-lined boundary represents the flying range of the drones when deployed from a satellite station. \label{fig:big-picture}}

    \resizebox{.5\textwidth}{!}{\begin{tikzpicture}[
roundnode/.style={circle, draw=blue!60, fill=green!2, very thick, minimum size=0.2mm},
roundnodeSmall/.style={circle, draw=red!60, fill=red!2, very thick, radius=0.05cm},
squarednode/.style={rectangle, draw=black!60, fill=black!2, very thick, minimum size=0.5mm},
->,>=stealth',shorten >=1pt,auto,node distance=1.5cm,semithick]

\node[roundnode] 	(newC2) at (4,2) {};
\node[roundnode] 	(newC4) at (1.3,-2.5) {};

 \draw[black,thick,dashed,fill=yellow!25!white,opacity=0.5] (newC2) circle (2.75);
\draw[black,thick,dashed,fill=yellow!25!white,opacity=0.5] (newC4) circle (2.75);
 \draw[black,thick,dashed,opacity=0.5] (newC2) circle (2.75);
\draw[black,thick,dashed,opacity=0.5] (newC4) circle (2.75);

\node[squarednode]        (depot)  at (-1,1.7)   {depot};
\node[roundnode] 	(C1) at (0.5,3) {1};
\node[roundnode, fill=blue!30!white] 	(C2) at (4,2) {2};
\node[roundnode] 	(C3) at (5,0) {3};
\node[roundnode, fill=blue!30!white] 	(C4) at (1.3,-2.5) {4};

\node[roundnodeSmall] 	(C9) at (3.25,3.25) {};
\node[roundnodeSmall, fill=red!30] 	(C5) at (4.25,3) {};
\node[roundnodeSmall, fill=red!30] 	(C6) at (5.4,2.2) {};
\node[roundnodeSmall, fill=red!30] 	(C7) at (4.6,4.0) {};
\node[roundnodeSmall, fill=green!50!black] 	(C8) at (3,4.25) {};
\node[roundnodeSmall] 	(Cx1) at (5,1.55) {};
\node[roundnodeSmall] 	(Cx11) at (6.25,1.75) {};
\node[roundnodeSmall,fill=green!50!black] 	(Cx2) at (2.3,3.65) {};
\node[roundnodeSmall] 	(Cx3) at (4.2,4.2) {};
\node[roundnodeSmall] 	(Cx4) at (3.6,4.5) {};

\node[roundnodeSmall, fill=green!50!black] 	(C13) at (2.5,-4) {};
\node[roundnodeSmall, fill=green!50!black] 	(C14) at (1.5,-4) {};
\node[roundnodeSmall] 	(Cx5) at (-0.1,-3) {};
\node[roundnodeSmall] 	(Cx6) at (0.15,-4.25) {};
\node[roundnodeSmall] 	(Cx7) at (-.7,-1.25) {};

\node[roundnodeSmall] 	(Cx17) at (3.25,-2.5) {};

\node[roundnodeSmall,fill=red!30] 	(C10) at (0.65,-3) {};
\node[roundnodeSmall,fill=red!30] 	(C11) at (0,-1) {};
\node[roundnodeSmall,fill=red!30] 	(C12) at (-1,-2.33) {};

\path[-] (depot) edge[black, thin] (C1);
\path[-] (C1) edge[black, thin] (C2);
\path[-] (C2) edge[black, thin] (C3);
\path[-] (C3) edge[black, thin] (C4);
\path[-] (depot) edge[black, thin] (C3);
\path[-] (C1) edge[black, thin] (C4);

\path (depot) edge[blue, very thick] (C2);
\path (C2) edge[blue, very thick] (C4);
\path (C4) edge[blue, very thick] (depot);

\path (C2) edge[red, dashed] (C5);
\path (C5) edge[red, dashed] (C7);
\path (C7) edge[red, dashed] (C6);
\path (C6) edge[red, dashed] (C2);

\path (C2) edge[green!50!black, dashed, bend left=15] (Cx2);
\path (Cx2) edge[green!50!black, dashed] (C8);
\path (C8) edge[green!50!black, dashed, bend left=15] (C2);

\path (C4) edge[green!50!black, dashed, bend right=0] (C13);
\path (C13) edge[green!50!black, dashed, bend right=0] (C14);
\path (C14) edge[green!50!black, dashed, bend right=0] (C4);

\path (C4) edge[red, dashed]                   (C11);
\path (C11) edge[red, dashed]                   (C12);
\path (C12) edge[red, dashed]                   (C10);
\path (C10) edge[red, dashed]                   (C4);

	\end{tikzpicture}}
\end{figure}
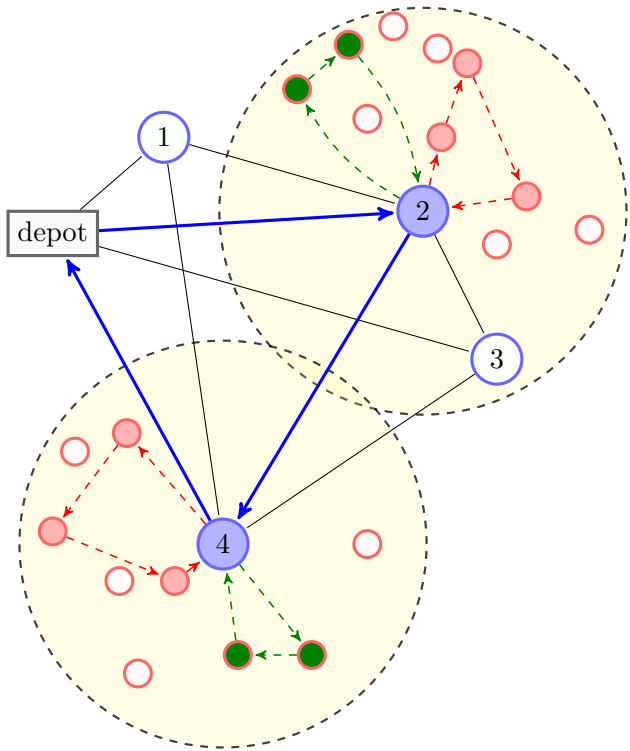

In the problem we aim to solve, the objective is to minimize the cost of delay in reaching communities and the cost of failing to reach a community or fully satisfy its (uncertain) demand, while addressing the following questions:

\begin{enumerate}
    \item Starting from the depot, what is the optimal set of satellite stations for a truck to visit?
    \item For each satellite station visited, what is the optimal set of affected communities for a delivery drone to serve?
\end{enumerate} We present a big-picture view of our problem in Figure \ref{fig:big-picture}. In this example, our one truck is asked to visit nodes 2 and 4 (through 3), and deploy two drones at each of the nodes to visit a subset of the affected communities before returning.

\begin{figure}
    \centering
    \caption{\label{fig:framework} The proposed two-stage robust decision framework.}
    \includegraphics[width=.75\textwidth]{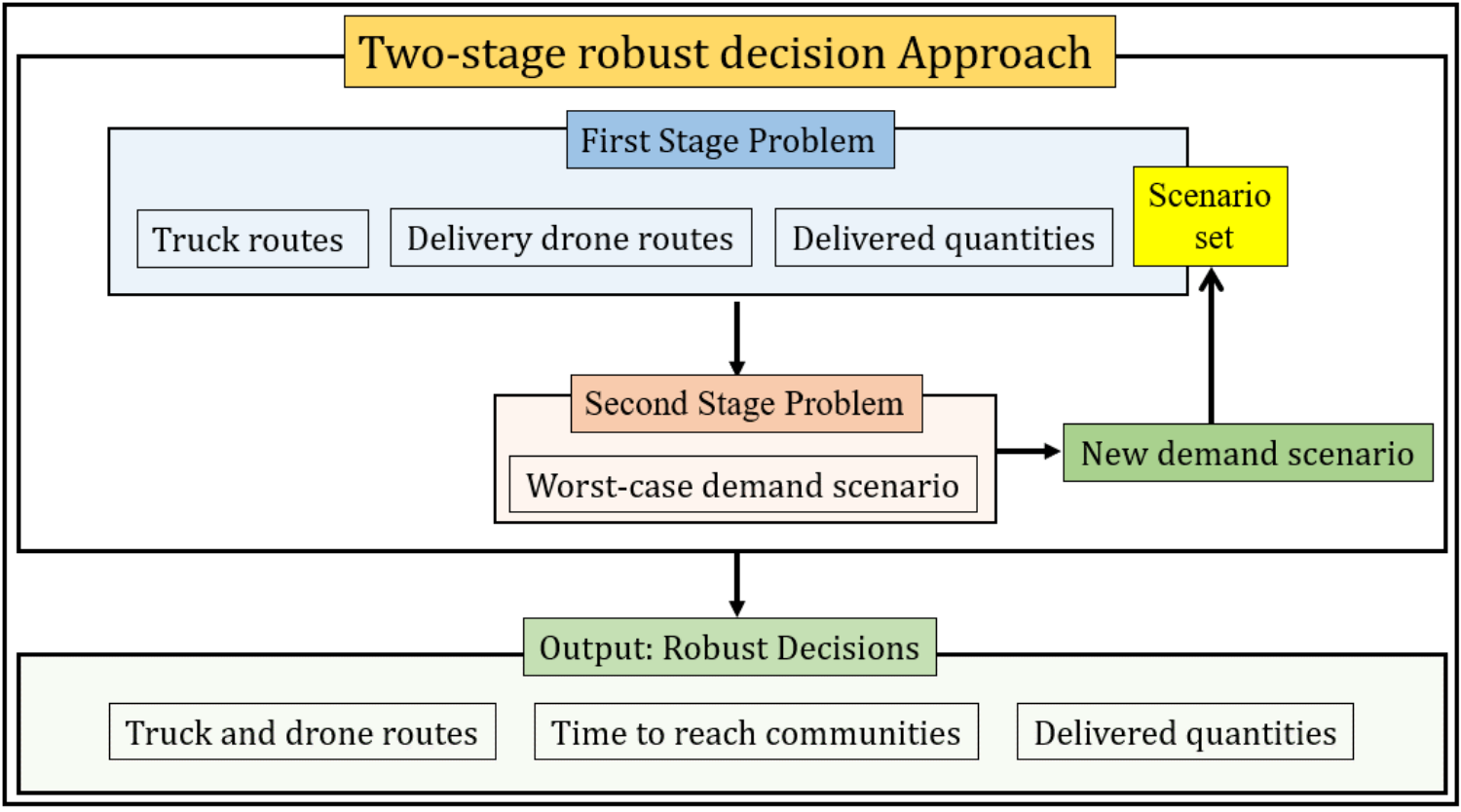}
\end{figure}

We now present the two-stage robust decision framework that we term R-2EVRP-UAV (\textbf{R}obust \textbf{2}-\textbf{E}chelon \textbf{V}ehicle \textbf{R}outing \textbf{P}roblem with \textbf{UAV}s). Contrary to other approaches in handling uncertainty (see, e.g., previous work where demands are first queried and then satisfied in \cite{our_paper}), here we propose a robust optimization framework. We first identify the worst-case demand instances given an uncertainty budget, and then make decisions for scheduling and routing to minimize the objective function. The framework is depicted for simplicity in Figure \ref{fig:framework}.

\begin{figure}[h!]
    \centering
        \caption{The proposed column generation for solving the first stage problem of R-2EVRP-UAV. In the main problem, we use a restricted set of the (exponentially many) drone routes, that we then supplement with new ones from our pricing subproblem.}
        \includegraphics[width=.9\textwidth]{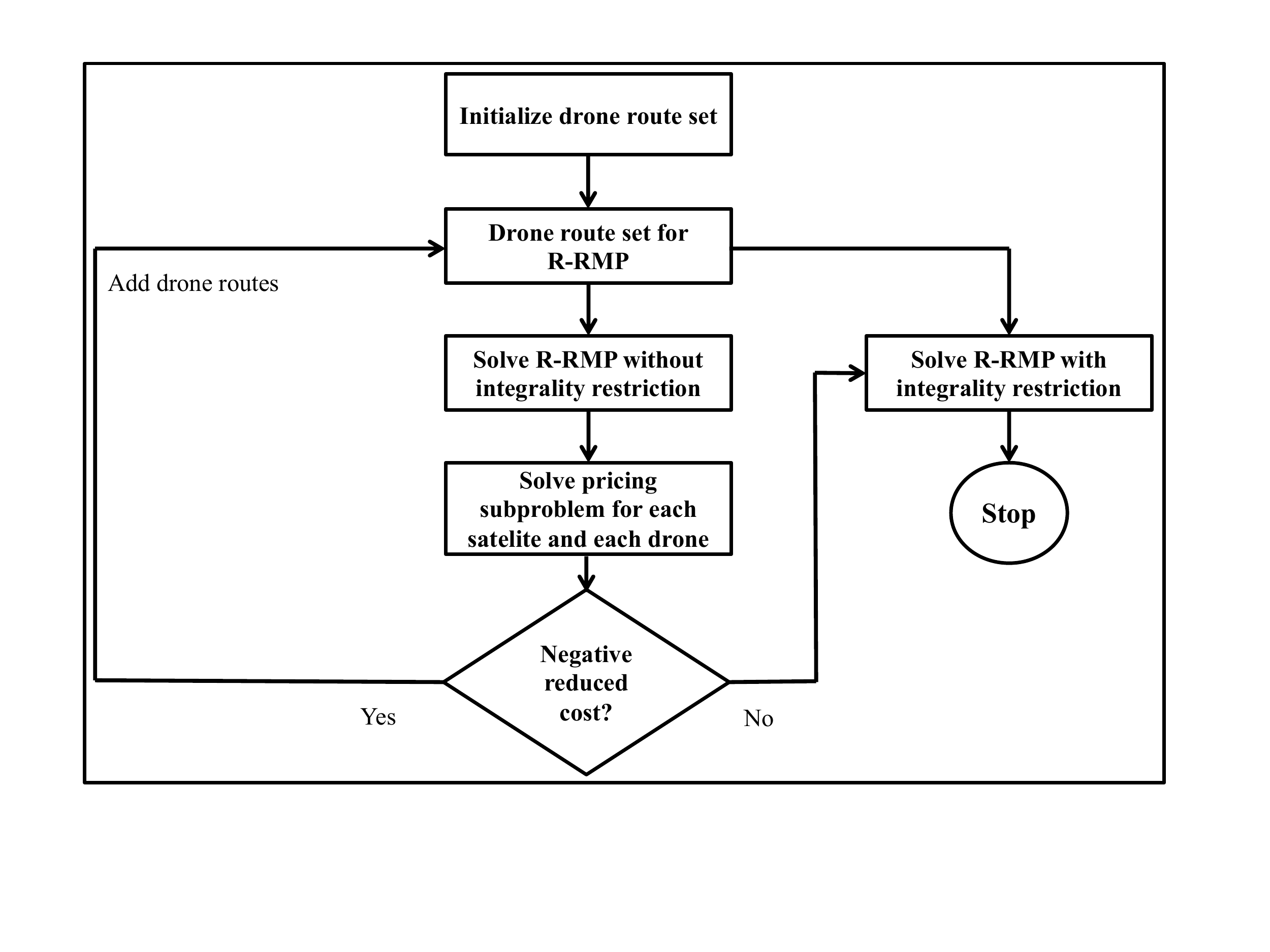}
    \label{fig:column_generation}
    \vspace{-0pt}
\end{figure}

The first stage problem is solved using a column generation framework. In the main problem, namely R-RMP (\textbf{R}obust \textbf{R}estricted \textbf{M}ain \textbf{P}roblem), we have a restricted set of drone routes that we can use, whereas in the pricing subproblem, we generate new routes to enter into our restricted set. We visualize this process, for convenience, in Figure \ref{fig:column_generation}.

\subsection{Robust Restricted main problem}
\label{subsec:main_problem}

We need the following decision variables. First, we have the first echelon (truck) routing variable as $x_{ij}$, which is equal to 1 if arc $\left(i,j\right)\in E^t$ is used. Second, we define a binary variable for the second echelon (drone/UAV) routing as $z_{p}^{s\ell}$, which is equal to 1 if drone $\ell\in D_s$, deployed from satellite station $s\in S$ is tasked with following route $p\in V_s$; 0, otherwise. As mentioned in the problem definition, some community $c\in C$ may be left without a delivery. If that happens, then binary variables $J_c$ will be equal to $1$. 

We also define timing variables. The elapsed time until we reach community $c\in C$ by some drone deployed from satellite station $s\in S$ is denoted as $T_{sc}^d$. With a slight misuse of notation, the elapsed time until a truck reaches a satellite station $s$ is denoted as $T_s$. The length of time that a truck spends idle at satellite station $s\in S$ is $\Delta_s$. 

Finally, we need variables to keep track of the cost for failing to deliver all quantities requested. This missed delivery amount at community $c\in C$ under scenario $\omega$ is defined as $R_c^\omega$. $\mathcal{Y}$ is an artificial variable representing the total recourse cost. With these definitions at hand, we now present the main problem in \eqref{main_problem}.

\allowdisplaybreaks

\begin{subequations}   
\label{main_problem}
\begin{align}
    \label{obj}
    \min~~ & \sum\limits_{s\in S} \sum\limits_{c\in C} F_c^T T_{sc}^d + \sum\limits_{c\in C} F_c^R J_c + \mathcal{Y} \\
    \label{number_of_trucks}
    \text{s.t.}~~ & \sum\limits_{j\in V:\left(0,j\right)\in E^t} x_{ij} \leq m^t, \\
    \label{flow_preservation}
    & \sum\limits_{j\in V:\left(i,j\right)\in E^t}x_{ij}-\sum\limits_{j\in V:\left(j,i\right)\in E^t}x_{ji}=0, & \forall i\in V\setminus\left\{0\right\}, \\
    \label{timing1}
    & T_j \geq T_i + \Delta_i + \tau_{ij}^t - M\left(1-x_{ij}\right), & \forall \left(i,j\right)\in E^t, \\
    \label{timing2}
    & T_j \leq T_i + \Delta_i + \tau_{ij}^t + M\left(1-x_{ij}\right), & \forall \left(i,j\right)\in E^t, \\
    \label{number_of_drones}
    & \sum\limits_{\ell \in D_s} \sum\limits_{p\in V_s} z_p^{s\ell}-m^d\sum\limits_{\left(i,s\right)\in E^t} x_{is}\leq 0, & \forall s\in S, \\
    \label{allowed_to_visit}
    & \sum\limits_{p\in V_s} z_p^{s\ell}-\sum\limits_{\left(i,s\right)\in E^t} x_{is} \leq 0, & \forall \ell\in D_s, \forall s\in S, \\
    \label{time_for_truck}
    & \sum\limits_{p\in V_s} \Pi_{sp}^d z_p^{s\ell} - \Delta_{s} \leq 0, & \forall \ell \in D_s, \forall s\in S, \\
    \label{visited_or_not}
    & \sum\limits_{s\in S}\sum\limits_{\ell \in D_s}\sum\limits_{p\in V_s} b_{c}^p z_p^{s\ell} + J_c = 1, & \forall c\in C, \\
    \label{time_for_community}
    & T_s+ \sum\limits_{\ell \in D_s} \sum\limits_{p\in V_s} \overline{t}_c^{s\ell}z_p^{s\ell}-M\left(1-\sum\limits_{\ell \in D_s}\sum\limits_{p\in V_s} b_{c}^p z_p^{s\ell}\right)\leq T_{sc}^d, & \forall c\in C_s, \forall s\in S, \\
    \label{recourse_constraint}
    & \mathcal{Y} \geq \sum\limits_{c\in C} F_c^D R_c^\omega, & \forall \omega \in \Omega^\prime, \\
    \label{missed_demand}
    & \sum\limits_{s\in S}\sum\limits_{\ell \in D_s}\sum\limits_{p\in V_s} z_p^{s\ell}d_{sc}^{\omega p}+R_c^\omega \geq Q_c^\omega, & \forall \omega \in \Omega^\prime, \forall c\in C, \\
    \label{variable_restriction_first}
    & x_{ij}\in\left\{0,1\right\}, & \forall \left(i,j\right)\in E^t, \\
    & z_p^{s\ell}\in\left\{0,1\right\}, & \forall \ell\in D_s, \forall p\in V_s, \forall s\in S, \\
    & J_c\in\left\{0,1\right\}, & \forall c\in C, \\
    & T_s\geq 0, & \forall s\in S, \\
    & \Delta_s\geq 0, & \forall s\in S, \\
    & T_{sc}^d\geq 0, & \forall c\in C_s, \forall s\in S,  \\
    & R_c^\omega \geq 0, & \forall c\in C, \forall \omega\in \Omega^\prime, \\
    \label{variable_restriction_last}
    & \mathcal{Y} \geq 0.
\end{align}
\end{subequations}

The objective function in \eqref{obj} minimizes: (a) the total time (delay cost) to reach all communities by delivery drones, (b) total cost of failure to reach communities, and (c) the recourse cost of failing to satisfy demands at the affected communities. Constraints \eqref{number_of_trucks} limit the maximum number of trucks that can be used; flow preservation for trucks is ensured by constraints \eqref{flow_preservation}. Constraints \eqref{timing1} and \eqref{timing2} dictate the time it takes a truck to reach each satellite node. In constraints \eqref{number_of_drones} and \eqref{allowed_to_visit}, we ensure that delivery drones are only dispatched from a satellite node if a truck has visited that node; furthermore, the number of delivery drones dispatched is limited by the number of drones carried by that truck. According to constraints \eqref{time_for_truck}, the lengths of time that a truck has to spend at a satellite location due to flights of delivery drones are bounded by the maximum flight times of delivery routes originating from that satellite, respectively. Constraints \eqref{visited_or_not} ensure that every community is either visited by a delivery drone or the variable $J_c$ takes a non-zero value. Constraint set \eqref{time_for_community} defines the required times to reach every community by a delivery drone from each satellite. Constraints \eqref{recourse_constraint} indicate that, for each scenario $\omega$, variable $\mathcal{Y}$ must take at least the value of the cost of failing to satisfy the demand quantities at that scenario. Constraints \eqref{missed_demand} ensure that, in each scenario, the demand at each affected community must be satisfied by the delivery drone or the recourse variable $R_c^\omega$ takes a positive value equal to the missed delivery amount. The variable restrictions are given in \eqref{variable_restriction_first}--\eqref{variable_restriction_last}, following the definitions provided. 

Finally, before we proceed to the pricing subproblem, we define dual variables associated with constraints \eqref{number_of_drones} ($\theta_s, \forall s\in S$),  \eqref{allowed_to_visit} ($\phi_{s\ell}, \forall \ell \in D_s, \forall s\in S$), \eqref{time_for_truck} ($\rho_{s\ell}, \forall \ell \in D_s, \forall s\in S$), \eqref{visited_or_not} ($\xi_c, \forall c\in C$), \eqref{time_for_community} ($\alpha_{sc}, \forall c\in C_s, \forall s\in S$), and \eqref{missed_demand} ($\beta_{c}^{\omega}, \forall \omega\in \Omega^\prime, \forall c\in C$). We can now discuss the corresponding pricing subproblem in the next subsection.

\subsection{Pricing subproblem}
\label{subsec:pricing}

As expected, we need some newly defined decision variables to formulate the subproblem. Seeing as in the pricing subproblem our goal is to generate routes (see Figure \ref{fig:column_generation}), we start with variable $v_{ij}^{s\ell}$, which is equal to $1$ if $\left(i,j\right)\in E^d_s$ is used by drone $\ell\in D_s$, and $0$ otherwise. The time it takes drone $\ell\in D_s$ before it reaches node $j\in C_s\cup \left\{s\right\}$ (i.e., a community that is reachable from the satellite station $s$ or the station itself at the end of the tour) is denoted by $t_j^{s\ell}$. The quantity transported on the drone $\ell\in D_s$ while traversing arc $\left(i,j\right)\in E_s^d$ under scenario $\omega$ is $q_{ij}^{\omega s \ell}$. Finally, the amount of items that is actually delivered at community $c\in C$ from drone $\ell\in D_s$ under scenario $\omega$ is $p_c^{\omega s\ell}$.

Let $RC_{s\ell}$ be the reduced cost for $s\in S, \ell\in D_s$. Using the dual variables obtained from \eqref{main_problem}, we have that: $$RC_{s\ell}=-\theta_s - \phi_{s\ell}-\sum\limits_{j\in C_s}\sum\limits_{\left(i,j\right)\in E_s^d}v_{ij}^{s\ell} \left(\xi_{sj}+\rho_{s\ell}\tau_{ij}^d\right)-\sum\limits_{j\in C_s}\alpha_{sj}\left(t_j^{s\ell}+M\sum\limits_{\left(i,j\right)\in E_s^d}v_{ij}^{s\ell}\right)-\sum\limits_{\omega \in \Omega^\prime}\sum\limits_{j\in C_s}\beta_{j}^{\omega}p_j^{\omega s \ell}.$$ Then, the pricing subproblem for every satellite station $s\in S$, and every drone $\ell \in D_s$ is given by the mixed integer linear program in \eqref{pricing_subproblem}. \\

\begin{subequations}
\label{pricing_subproblem}
    \begin{align}
        \label{pricing_obj}
        \min~~ & RC_{s\ell} \\
        \label{pricing_range}
        \text{s.t.}~~ & \sum\limits_{\left(i,j\right)\in E_s^d} v_{ij}^{s\ell} \tau_{ij}^d \leq W^d, \\
        \label{pricing_at_most_one_route}
        &\sum\limits_{\left(i,j\right)\in E_s^d} v_{sj}^{s\ell} \leq 1, \\
        \label{pricing_at_most_one_community}
        &\sum\limits_{\left(i,j\right)\in E_s^d} v_{ij}^{s\ell} \leq 1, & \forall j\in C_s, \\
        \label{pricing_flow_preservation}
        & \sum\limits_{i:\left(i,j\right)\in E_{s}^d} v_{ij}^{s \ell}-\sum\limits_{i:\left(j,i\right)\in E_{s}^d} v_{ji}^{s \ell}=0, & \forall j\in C_s, \\
        \label{pricing_timing1}
        & t_{j}^{s\ell} \geq t_{i}^{s\ell}+\tau_{ij}^d-W^d\left(1-v_{ij}^{s\ell}\right), & \forall \left(i,j\right)\in E_s^d, \\
        \label{pricing_timing2}
        & t_{j}^{s\ell} \leq t_{i}^{s\ell}+\tau_{ij}^d+W^d\left(1-v_{ij}^{s\ell}\right), & \forall \left(i,j\right)\in E_s^d, \\
        \label{pricing_max_quantity1}
        &q_{ij}^{\omega s \ell} \leq v_{ij}^{s \ell} L_{max}, & \forall \omega \in \Omega^\prime, \forall \left(i,j\right)\in E_s^d, \\
        \label{pricing_quantity_carried}
        &\sum\limits_{\left(i,j\right)\in E_s^d} q_{ij}^{\omega s \ell}-\sum\limits_{\left(j,k\right)\in E_{s}^d}q_{jk}^{\omega s \ell}=p_{j}^{\omega s \ell}, & \forall \omega \in \Omega^\prime, \forall j \in C_s, \\
        \label{pricing_max_quantity2}
        &\sum\limits_{j\in C_s} p_j^{\omega s \ell}\leq L_{max}, & \forall \omega \in \Omega^\prime, \\
        \label{pricing_variable_restrictions1}
        &v_{ij}^{s\ell}\in \left\{0,1\right\}, & \forall \left(i,j\right)\in E_s^d, \forall \ell \in D_s, \\
        & t_j^{s \ell} \geq 0, & \forall j\in C_s, \\
        & q_{ij}^{\omega s \ell}\geq 0, & \forall \omega \in \Omega^\prime, \forall \left(i,j\right)\in E_s^d, \\
        \label{pricing_variable_restrictions2}
        & p_{j}^{\omega s \ell}\geq 0, & \forall \omega \in \Omega^\prime, \forall j\in C_s.
    \end{align}
\end{subequations}

The objective function in \eqref{pricing_obj} minimizes the reduced cost of a delivery drone route. That way, we can pick new routes to add to the main problem (R-RMP). All new routes need to respect the flying range capabilities of a drone, and this is enforced in \eqref{pricing_range}. Constraints \eqref{pricing_at_most_one_route} and \eqref{pricing_at_most_one_community} state that a valid route is one that starts from the satellite station and visits each community location at most once. Constraints \eqref{pricing_flow_preservation} are traditional flow preservation constraints at each of the community locations. Constraints \eqref{pricing_timing1} and \eqref{pricing_timing2} define the time that a drone can reach a community. Constraints \eqref{pricing_max_quantity1}--\eqref{pricing_max_quantity2} work together to ensure that the drone does not carry more load than its capacity $L_{max}$ and that the amount of product that a drone carries is reduced by the amount  delivered in communities visited earlier in the route.

We note here that the subproblem is posed for every satellite station $s\in S$ and every UAV $\ell\in D_s$ that originates from a truck stopped at $s$. This implies that the best route for the first delivery drone picked from $D_s$ would also be the best route for every other drone in $D_s$. To avoid generating the same route multiple times, and in order to identify a diverse set of candidate columns to enter into the restricted main problem, we design a heuristic approach by updating the list of communities to be visited for each subsequent pricing subproblem. 

Let $\mathcal{C}^s_\ell$ be the set of community nodes visited when solving the pricing subproblem for satellite station $s$ and drone $\ell$. In mathematical terms, we have $$\mathcal{C}^s_\ell=\left\{j\in C_s:\sum\limits_{\left(i,j\right)\in E_s^d} v_{ij}^{s\ell}=1\right\}.$$ The update is quite simple: we drop all communities that are included in previously generated routes from $C_s$. Mathematically, for pricing subproblem $\ell$, we update as follows: \begin{align}\label{update_community_set}C_s \leftarrow C_s\setminus \left\{\bigcup\limits_{\ell^\prime <\ell}C_{s}^{\ell^\prime}\right\},\end{align} after having solved all previous corresponding pricing subproblem for UAVs $\ell^\prime<\ell$.

As an example, for the first delivery drone in $D_s$ we will consider the best route that can pass through any community in $C_s$. For the second delivery drone, though, we heuristically update $C_s$ as $C_s\setminus C_{s}^1$, and then we further update $C_s$ to $C_s\setminus \left(C_s^1\cup C_s^2\right)$, and so on.

This results in obtaining solutions to the pricing subproblems faster, as well as obtaining distinct routes to use as columns in R-RMP, at the expense of higher route cost estimations.

\subsection{Worst-case demand scenario generation}
\label{subsec:scenario_generation}

In this subsection, we discuss the scenario generation subproblem, which determines the worst-case demand scenario for the current solution of the restricted main problem (R-RMP), the demand distribution, and the uncertainty budget. 

For every community $c\in C$, let $\overline{Q}_c$ be the average demand and let $\hat{Q}_c$ be the maximum demand deviation that can be experienced in that location. The average demand can be estimated from census data that includes the number of people in a community, their socioeconomic information, their age. On the other hand, $\hat{Q}_c$ is dictated by the origin of stochasticity in the problem. For example, disaster type, severity, or movement direction/progression will affect the deviations we expect in post-disaster aid demand. We define $y_c$ as a binary variable that indicates whether a community $c$ experiences this maximum deviation in demand. We also define $\Gamma$ as the total uncertainty budget, as in the number of communities that can experience maximum deviations in demand simultaneously. 

Finally, we assume that the affected area can be divided into a set of geographical regions $A$. We group the communities that fall in a region (i.e., are located within their boundary) in $C^a$, such that $C=\bigcup\limits_{a\in A}C^a$ and $C^{a}\cap C^{a^\prime}=\emptyset$ for $a\neq a^\prime$. Now, we similarly define $\Gamma_a$ for every geographical region $a\in A$ as the uncertainty budget per geographical region. Introducing geographical regions, their associated boundaries, and individual uncertainty budgets $\Gamma_a$ for them, is necessary in the disaster management setting, as the impact of a disaster can vary significantly depending on the disaster characteristics and the location of a community. For example, depending on the direction and speed of movement of a hurricane, some areas may face the full-force of the devastation, while others may remain relatively safe. 


The uncertainty set can be written as: $$\mathcal{U}_Q=\left\{ Q\in\mathbb{R}^{|C|}: \Tilde{Q}_c=\overline{Q}_c+y_c\Hat{Q}_c, y_c\in\left\{0,1\right\} \forall c\in C, \sum\limits_{c\in C}y_c\leq \Gamma, \sum\limits_{c\in C^a} y_c\leq \Gamma_a, \forall a\in A \right\}.$$ Since our problem has simple recourse, the second stage problem for determining the most violated scenarios can be formulated as a maximization problem. The optimal solution of the problem is denoted by $OPT(\Gamma_Q)$.

Letting $\Tilde{Q}_c=\overline{Q}_c+y_c\Hat{Q}_c$ and introducing $\Psi_{c}^{s\ell}=\max\limits_{\omega \in \Omega^\prime}\sum\limits_{p\in V_s}z_{p}^{s\ell}d_{sc}^{\omega p}$ (i.e., the maximum amount delivered at community $c\in C$ by a UAV $\ell\in D_s$ following a route $p\in V_s$ starting from satellite $s\in S$ over all current scenarios $\omega \in \Omega^\prime$), we have the following max-min problem for generating the worst-case scenario for the current solution: \begin{subequations}\begin{align}
OPT(\Gamma_Q)=&\max\limits_{\Tilde{Q}\in\mathcal{U}_Q}~~  \min\limits_{R_c} \sum\limits_{c\in C} F_c^D R_c \\
\label{dualize_this}
& \text{s.t.}~~~~ R_c+\sum\limits_{s\in S}\sum\limits_{\ell \in D_s} \Psi_c^{s\ell} \geq \overline{Q}_c+\Hat{Q}_c, & \forall c\in C, \\
& ~~~~~~~~ R_c \geq 0, & \forall c\in C.
\end{align}
\end{subequations}

We observe that this bilevel formulation can be rewritten as a single maximization problem by dualizing the inner minimization problem. Let $\pi_c$ be the dual variables corresponding to constraints \eqref{dualize_this}. By the definition of the uncertainty set $\mathcal{U}_Q$ we have that $\Tilde{Q}_c=\overline{Q}_c+y_c\Hat{Q}_c$, which results in non-linearity due to the quantity $y_c\pi_c$ in the dualized inner problem. Hence, we introduce $\delta_c=y_c\pi_c$ and the corresponding linearization constraints (see \eqref{linearization1}--\eqref{linearization4}) to obtain the final mixed integer linear program of \eqref{final_mip}:
\begin{subequations}
\label{final_mip}
    \begin{align}
         \max~~& \sum\limits_{c\in C} \left(\overline{Q}_c\pi_c+\delta_c \Hat{Q}_c\right)-\sum\limits_{s\in S}\sum\limits_{\ell \in D_s}\sum\limits_{c\in C}\pi_c\Psi_c^{s\ell} \\
         \text{s.t.}~~& \sum\limits_{c\in C^a} y_c \leq \Gamma^a, & \forall a\in A, \\
         & \sum\limits_{c\in C} y_c \leq \Gamma, \\
         \label{linearization1}
         & \pi_c \leq F_c^D, & \forall c\in C, \\
         \label{linearization2}
         & \delta_c \leq \pi_c, & \forall c\in C, \\
         \label{linearization3}
         & \delta_c \leq F_c^D y_c, & \forall c\in C, \\
         \label{linearization4}
         & \delta_c \geq \pi_c - F_c^D\left(1-y_c\right), & \forall c\in C, \\
         &y_c\in\left\{0,1\right\}, \pi_c \geq 0, \delta_c \geq 0, & \forall c\in C.
    \end{align}
\end{subequations}

\subsection{Algorithm}
\label{sec:algorithm}

\begin{figure}[htbp]
    \centering
    \caption{A schematic representation of the algorithm used to solve the robust 2EVRP problem we are tackling in our framework.}
    \label{fig:schematic}
    \includegraphics[width=.99\textwidth]{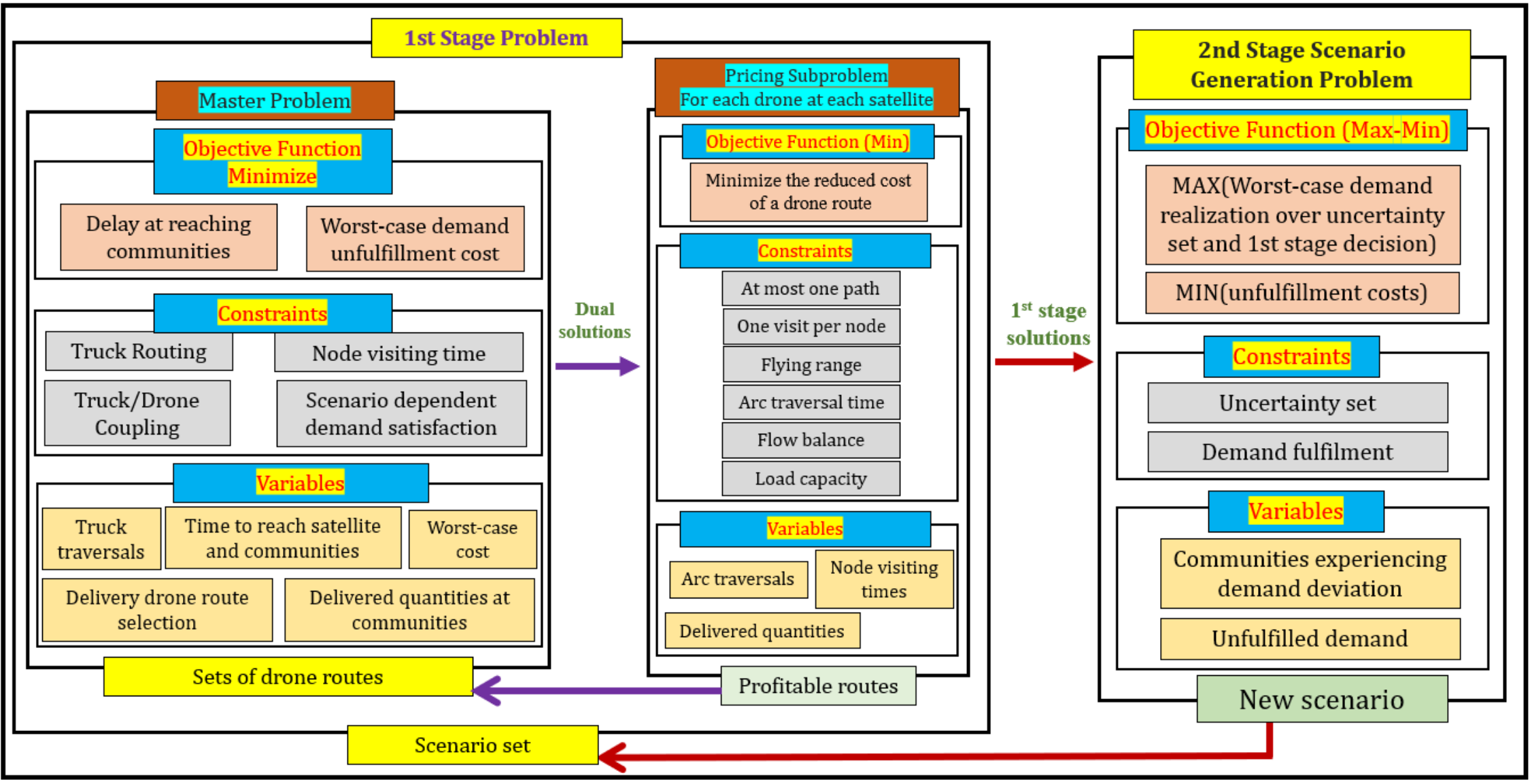}
\end{figure}

We finish this subsection by recounting the steps on how the framework is solved. We also present a schematic flow diagram of the algorithmic steps in Figure \ref{fig:schematic}.

\paragraph*{Step 1} We initialize an iteration counter ($n$), a lower ($LB$) and an upper bound ($UB$), and provide an initial scenario set ($\Omega^\prime$). Specifically, we set $n=1, LB=0, UB=+\infty$, and $\Omega^\prime=\Hat{\Omega}$, where $\Hat{\Omega}$ represents a scenario where each community has its average demand level $\overline{Q}_c$. We also generate an initial set of feasible routes for the delivery drones and add them to $V_s, \forall s\in S$.

\paragraph*{Step 2} We relax the variable integrality restrictions and solve a relaxed version of the restricted main problem (R-RMP from subsection \ref{subsec:main_problem}) with the current set of routes $V_s$ and the current set of scenarios $\Omega^\prime$. Solving, we obtain the duals $\theta_s$, $\phi_{s\ell}$, $\xi_c$, $\rho_{s\ell}$, $\alpha_{sc}$, $\beta_c^\omega$.

\paragraph*{Step 3} We solve the pricing subproblem (from subsection \ref{subsec:pricing}) for each satellite $s\in S$ and each delivery drone $\ell\in D_s$ using the current scenario set $\Omega^\prime$ and the most recently obtained dual variable values. As described in subsection \ref{subsec:pricing}, we solve the subproblem for each drone in set $D_s$ and update $C_s$ as in \eqref{update_community_set} until either $C_s=\emptyset$ or we have generated $m^d$ routes. If $\exists s,\ell$ such that $RC_{s\ell}\geq 0$ (i.e., with desired reduced cost), we add the corresponding routes to $V_s$ in R-RMP.

\paragraph*{Step 4} If $RC_{s\ell}\geq 0$ for every $s\in S, \ell\in D_s$, then we go to Step 5. Otherwise, assuming there exists some $\hat{s}, \hat{\ell}$ such that $RC_{\hat{s}\hat{\ell}}<0,$ we go back to Step 2 and resolve with the updated restricted set of routes $V_s$.

\paragraph*{Step 5} We solve the relaxed R-RMP again, but now with the updated restricted route set $V_s$ and the current scenario set $\Omega^\prime$. The main problem is now solved using the relax-and-fix heuristic. First, we solve by relaxing the variables for the delivery drones ($z$ space) and obtain the binary solution for the truck routes ($x$ space). Next, we fix the solutions for $x$ and we solve with the integrality of $z$ restored. The solutions to $z$ are used to update $\Psi_{c}^{s\ell}$ as $\Psi_{c}^{s\ell}=\max\limits_{\omega \in \Omega^\prime}\sum\limits_{p\in V_s}z_{p}^{s\ell}d_{sc}^{\omega p}$. The optimal objective function value for the restricted main problem is used as our new lower bound now, i.e., $$LB=\sum\limits_{s\in S} \sum\limits_{c\in C} F_c^T T_{sc}^d + \sum\limits_{c\in C} F_c^R J_c + \mathcal{Y}.$$

\paragraph*{Step 6} We move to the second stage and solve the scenario generation subproblem. The obtained optimal objective function value $OPT(\Gamma_Q)$ is used to update the upper bound as follows: $$UB\leftarrow \min\left\{UB, ~ \sum\limits_{s\in S} \sum\limits_{c\in C} F_c^T T_{sc}^d + \sum\limits_{c\in C} F_c^R J_c+OPT(\Gamma_Q)\right\}.$$

\paragraph*{Step 7} We now check our termination criterion based on the current values for our bounds. Specifically, we separate our discussion in two parts, depending on $UB-LB$. 

\vspace{-15pt}

\paragraph*{\indent \indent Step 7a} For the current bounds, if $UB-LB > \epsilon$, where $\epsilon$ is a user-defined tolerance parameter, then the new demand scenario $\omega$ is added to $\Omega^\prime$. The corresponding demand of this scenario is given by $Q_c^\omega=\overline{Q}_c+y_c \Hat{Q}_c$. A new column is added for the variables $R_c^\omega$ and new rows are added according to constraints \eqref{recourse_constraint}--\eqref{missed_demand}, and we go back to Step 2.

\vspace{-15pt}

\paragraph*{\indent \indent Step 7b} If $UB-LB\leq \epsilon$, then the solutions with the current scenario set $\Omega^\prime$ are the final robust truck and drone scheduling and routing decisions that can cover the worst-case demand scenario for the given demand distribution and uncertainty budget at minimum cost.

\section{Results}
\label{sec:results}

As mentioned in the Introduction, we are motivated by the humanitarian relief operations that need to take place in an area with numerous failures in its transportation, power, and communication networks. Hence, we present the applicability of our framework with a small case study based on generated data simulating emergency aid demand in the island of Puerto Rico after Hurricane Maria. In the following sections, we discuss our data generation process; we present our findings on computational runtime and measures of success of the operations, and perform a sensitivity analysis on drone parameters, i.e., the number of drones per truck and the drone flying range capability.

\subsection{Data generation}

We use the same datasets and data generation process as in \cite{our_paper}; that said, we introduce new parameters that are necessary to evaluate the robust framework. Specifically, we solve problems with number of communities $|C|\in \left\{60,80,100,120\right\}$. Each community represents a zip code of Puerto Rico. Since there are 173 zip codes in total, for each different problem size (i.e., different value of $|C|$) we randomly select out of those 173 zip codes. After selecting the communities that will have failed communications, we obtain average demand quantities $\overline{Q}_c$ by normalizing the zip code population into a range of $\left[2,15\right]$ units. Finally, the maximum demand deviation that community $c$ can experience, $\hat{Q}_c$ is set at $0.5\overline{Q}_c$, i.e., 50\% of the expected demand.

For the truck and drone details, we assume that $m^t=4$ (that is, 4 trucks), but we let $m^d$ be equal to 6, 9, 10, and 14 drones per truck for the different problem sizes (60, 80, 100, and 120 affected communities, respectively). $L_{max}$, the maximum allowed carry load is set at 25 units; $W^d$, the flying range, is set at 35 miles. That said, we do investigate the effect of the number of drones per truck in subsection \ref{subsec:number_of_drones} as well as the range $W^d$ in subsection \ref{subsec:range_of_drones}. All distances that the drones traverse are calculated based on the corresponding endpoints geodesic distances (``as the crow flies''). Both truck and drone velocities are deterministic and equal to 60 miles per hour. We note here that drone velocities are not deterministic, nor constant, and rely on the weather conditions among other environmental characteristics. Velocities and ranges have been identified as important future research avenues in the literature \cite{poikonen2021future} and in our future work recommendations in Section \ref{sec:conclusions}.

We then discuss cost selection. We have three costs associated with our operations: (a) delivery delay cost $F_c^T$ which is set at $\$1/minute$, (b) failure to reach a community $F_c^R$ which is set at $\$10000/community$, and (c) failure to fulfill all of the demand in a community $F_c^D$ which is randomly generated in $\$\left[10,1000\right]$ per unit for each community. This aims to address that some items are very important (e.g., epi-pens) and failing to fulfill them will have more serious repercussions. 

We finish the data description with the uncertainty budget. We consider $\Gamma\in\left\{30\%,50\%,70\%\right\}$. For the individual geographical region uncertainty budget $\Gamma^a, \forall a\in A$, we divide Puerto Rico into 10 mutually exclusive subregions based on the longitudinal values (from east to west) and then set $\Gamma^a=50\%$ for each subregion.

\subsection{Numerical experiments}

We implemented the framework using AMPL and CPLEX. We use a Dual Intel Xeon Processor (12 Core, 2.3GHz Turbo) computer with 64 GB of RAM. Furthermore, for each community size ($\lvert C\rvert$) and budget uncertainty ($\Gamma$), we generate 5 instances at random. We present the mean and standard deviation of our results from all of these experiments in Table \ref{tab:results}. We present five interesting metrics for the performance of our framework:
\begin{enumerate}
    \item CPU Time (in seconds) reveals the computational time required.
    \item \# Scenarios presents the number of scenarios that are generated to identify the optimal robust decisions.
    \item Cost (in \$) is the objective function value obtained.
    \item Unfilfilled demand (in \%) keeps track of the demand that was unfulfilled after performing the robust routes.
    \item Delay (in minutes) shows the time it takes until a community is reached.
\end{enumerate} For all the metrics, we present both the average from $5$ experiments as well as the standard deviation.

\newcolumntype{P}[1]{>{\centering\arraybackslash}p{#1}}

\begin{table}[htbp]
    \centering
    \caption{\label{tab:results} A summary of the results from the numerical experiment for different number of affected communities ($\lvert C \rvert$) and uncertainty budget ($\Gamma$). For all results, we show the average and standard deviations from 5 different runs.}
    \resizebox{\textwidth}{!}{\begin{tabular}{c|c|ll|ll|ll|ll|ll}
    \toprule
         \multicolumn{2}{c|}{\textbf{Instance}}  & \multicolumn{2}{c|}{\textbf{CPU time~(sec)}} & \multicolumn{2}{c|}{\textbf{\# Scenarios}} & \multicolumn{2}{c|}{\textbf{Cost (\$)}} & \multicolumn{2}{P{0.15\textwidth}|}{\textbf{Unfulfilled demand (\%)}} & \multicolumn{2}{c}{\textbf{Delay (min)}} \\
         \multicolumn{1}{c}{$\lvert C \rvert$}  & \multicolumn{1}{c|}{$\Gamma$ (\%)} & \multicolumn{1}{c}{$\mu$} & \multicolumn{1}{c|}{$\sigma$} & \multicolumn{1}{c}{$\mu$} & \multicolumn{1}{c|}{$\sigma$} & \multicolumn{1}{c}{$\mu$} & \multicolumn{1}{c|}{$\sigma$} & \multicolumn{1}{c}{$\mu$} & \multicolumn{1}{c|}{$\sigma$} & \multicolumn{1}{c}{$\mu$} & \multicolumn{1}{c}{$\sigma$}  \\
         \hline 
         \multirow{3}{*}{60} & 30 & 1131.67 &	786.37 &	4.33 &	0.47 &	11112.17 &	6355.38	& 5.15 &	5.42 &	63.33 &	8.87 \\
 & 50 & 784.33 &	515.94 &	3.50 &	0.50 &	8919.17 &	3538.31 &	3.65 &	3.38 &	57.75 &	11.50 \\
& 70 & 526.50 &	325.08 &	3.17 &	0.37 &	8467.50 &	5226.64 &	3.30 &	4.36 &	60.04 &	13.93 \\
\midrule 
\multirow{3}{*}{80} & 30 & 2504.83	 &789.42 &	4.17 &	0.37 &	11571.50 &	2545.54 &	4.61 &	2.01 &	71.56 &	7.58 \\
 & 50 & 1024.17 &	320.07 &	3.00 &	0.00 &	9226.17 &	2982.80 &	2.56 &	2.53 &	72.58 &	8.17 \\
& 70 & 1141.00 &	320.88 &	3.17 &	0.37 &	9232.00 &	2412.41 &	2.59 &	1.89 &	73.55	 &9.03 \\
\midrule 
\multirow{3}{*}{100} & 30 & 4050.00 &	1643.08 &	4.33 &	0.47 &	12252.00 &	1320.08 &	3.12 &	0.75 &	100.04	 &11.69 \\
& 50 & 1812.33 &	581.74 &	3.00 &	0.00 &	9649.67 &	1739.32 &	1.42 &	0.91 &	92.08 &	9.98 \\
& 70 & 1773.67 &	531.38 &	3.17	 &0.37	 &9490.00 &	1838.55 &	1.38 &	0.98 &	90.88	 &8.61 \\
\midrule 
\multirow{3}{*}{120} & 30 & 5417.33	 &1580.40	 &4.17	 &0.69	 &16341.50 &	954.90	 &3.46	 &1.45	 &131.52 &	34.06 \\
& 50 & 3712.33 &	1362.32 &	3.17 &	0.37	 &11204.83	 &1934.51 &	1.33 &	0.98	 &120.13 &	19.60 \\ 
& 70 & 3907.83 &	1615.97 &	3.33	 &0.75	 &11561.67	 &2780.54 &	1.43	 &1.18 &	120.74 &	29.48 \\
\bottomrule
    \end{tabular}}
\end{table}

We see that in all of our experiments, the number of scenarios required before we can identify the optimal robust decisions is very small, both on average and when considering the deviations. We also observe that problem instance size does not really play a role in the number of scenarios generated. 

In addition, all three metrics (cost of operations, proportion of unfulfilled demand, delay in reaching a community) all decrease as the uncertainty budget increases. Here, average proportion of unfulfilled demand is calculated by taking the average of unfulfilled demand proportions in all considered scenarios. More specifically, we find that an uncertainty budget of 30\% results in higher values for all these metrics. When the uncertainty budget is 50\% and 70\%, the difference between the average values for total cost, proportion of unfulfilled demand, and average delay time is very small. This is attributed to the area uncertainty budget $\Gamma^a$, which is set at 50\%.

\subsection{Sensitivity analysis}

We also design experiments to investigate the effect of drone parameters on the routing decisions and the metrics identified earlier (operational cost, proportion of unfulfilled demand, delay to community). Specifically, we want to investigate whether more drones per truck would result in faster, more efficient operations. We also investigate whether the flying range of the drones is an important parameter. 

\begin{figure}
    \centering
    \caption{The four (geographically organized) randomly generated instances for the: (a) northwestern part, (b) northeastern part, (c) southwestern part, and (d) southeastern part of Puerto Rico. \label{fig:datasets}}
    \includegraphics[width=.99\textwidth]{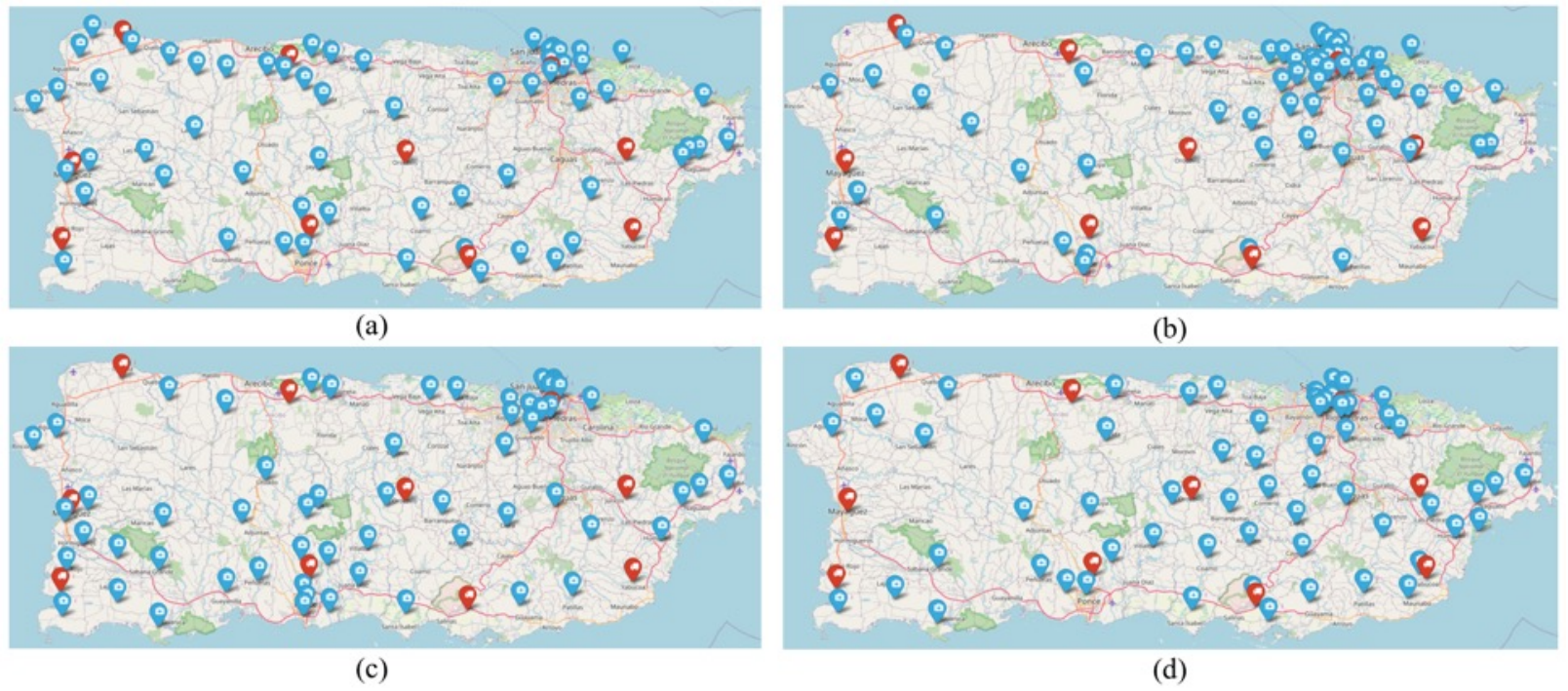}
\end{figure}

To do that, we focus on $\lvert C \rvert=60$ communities per instance and generate 5 different instances, per geographical region. The first four instances are randomly generated by selecting communities across Puerto Rico, however with higher probability in specific areas at a time. The four instances focus on the: (a) northwestern part of the island, (b) northeastern part of the island, (c) southwestern part of the island, and (d) southeastern part of the island, respectively. The instances generated for this experiment are presented in Figure \ref{fig:datasets}. The fifth instance is constructed by randomly selecting locations without any geographical preference (i.e., all locations are equally probably selected).

We also introduce a new attribute, i.e., disaster level to represent the disaster impact/strength, which in turn controls the maximum demand deviation a geographical area can experience. For simplicity, we only consider three levels of disaster impact. In the first level (Level 1), communities can experience maximum deviation that is 50\% of the average demand irrespective of their geographical locations. In Levels 2 and 3, the geographical location of a community dictates the maximum deviation that it can experience. We use the same 10 geographical areas introduced earlier (i.e., Puerto Rico consists of 10 mutually exclusive subregions based on the longitudinal values). In disaster level 2, the highest demand deviation (90\% of the average demand) is experienced by the eastern-most area, and it gradually decreases (by 10\% from one area to the next) from east to west, with the least deviation (0\% of the nominal demand) being experienced by the western-most area. In disaster level 3, we reverse the direction. Then, we set the area uncertainty budget $\Gamma^a$ at 100\%, which indicates that every community in geographical area $a$ can experience simultaneous demand deviation in a scenario. 

\subsubsection{Impact of number of drones per truck} \label{subsec:number_of_drones}

We test three cases for the number of drones per truck: $m^d=4$, $m^d=6$, and $m^d=8$. The main result is that (as expected), as the number of drones per truck increase, all metrics improve. That said, it is interesting to note that the rate of improvement depends on the disaster impact level and the uncertainty budget. 

\begin{figure}
    \centering
    \caption{Variation in average unfulfilled demand and average delay to reach a community as the number of delivery drones per truck changes. \label{fig:number_of_drones}}
    \includegraphics[width=.99\textwidth]{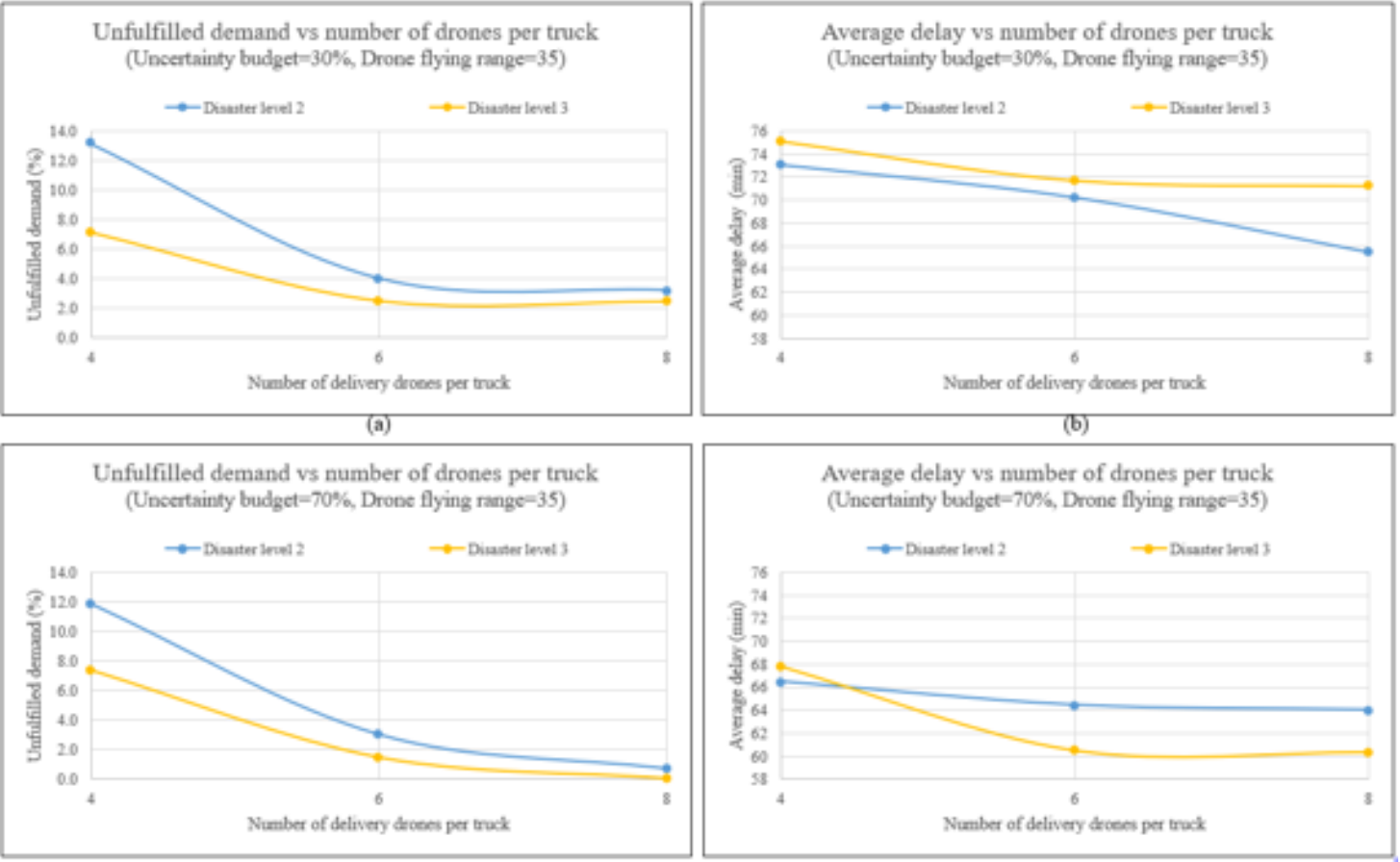}
\end{figure}

Specifically, in Figure \ref{fig:number_of_drones}, we observe that the rate of decrease in the average unfulfilled demand with the number of delivery drones is very similar for Levels 2 and 3. However, for disaster level 2, a higher proportion of demand is unfulfilled. This is due to geography of Puerto Rico, where a bigger population lives in communities in the eastern areas. As a reminder, disaster level 2 results in higher demand deviations in these areas. 

We also observe a slight difference in the rate of decrease in unfulfilled demand between uncertainty budget of 30\% and 70\% in Figures \ref{fig:number_of_drones}(a) and \ref{fig:number_of_drones}(c). When the uncertainty budget is equal to 30\%, we see a decrease in the average proportion of unfulfilled demand when increasing the number of drones from 4 to 6. That said, adding more drones per truck only slightly reduce the unfulfilled demand. We discussed earlier that lower uncertainty budget results in higher number of scenarios to be generated, where each scenario is characterized by fewer communities experiencing demand deviation. This results in demand scenarios, where fulfilment is more constrained by the drone flying range rather than by the drone load capacity. Without increasing the flying range, additional drones only slightly improve the demand fulfilment. 

However, we see in Figure \ref{fig:number_of_drones}(b) that there is a reduction in the average delay in reaching communities as we increase the number of drones. Again, we see different rates of decrease in the average delay: for disaster level 3, where communities in the western areas experience higher demand deviation, the rate of decrease is very small. This is because the communities in the western areas are sparsely located and farther apart. Hence, they require longer times to reach from nearby satellite stations. This is why increasing the number of drones does not necessarily result in a reduction in the average delay time. In contrast, for disaster level 2, more densely packed communities experience higher demand deviation and with the increase of number of drones comes a large reduction in the average delay to reach a community. For the case of uncertainty budget equal to 70\%, few scenarios with communities experiencing high demand deviation are generated. The resulting demand scenarios are more constrained by drone load capacity than by flying range. As the number of drones is increased from 4 to 6, the proportion of unfulfilled demand and the average delay time both decreases. The proportion of unfulfilled demand is close to 0\% when more drones are added while the corresponding average delay time remains almost unchanged.

\subsubsection{Impact of range of drones} \label{subsec:range_of_drones}

\begin{figure}[htbp]
    \centering
    \caption{Variation in average unfulfilled demand and average delay with delivery drone flying range (8 delivery drones per truck). \label{fig:range_of_drones_8}}
    \includegraphics[width=.99\textwidth]{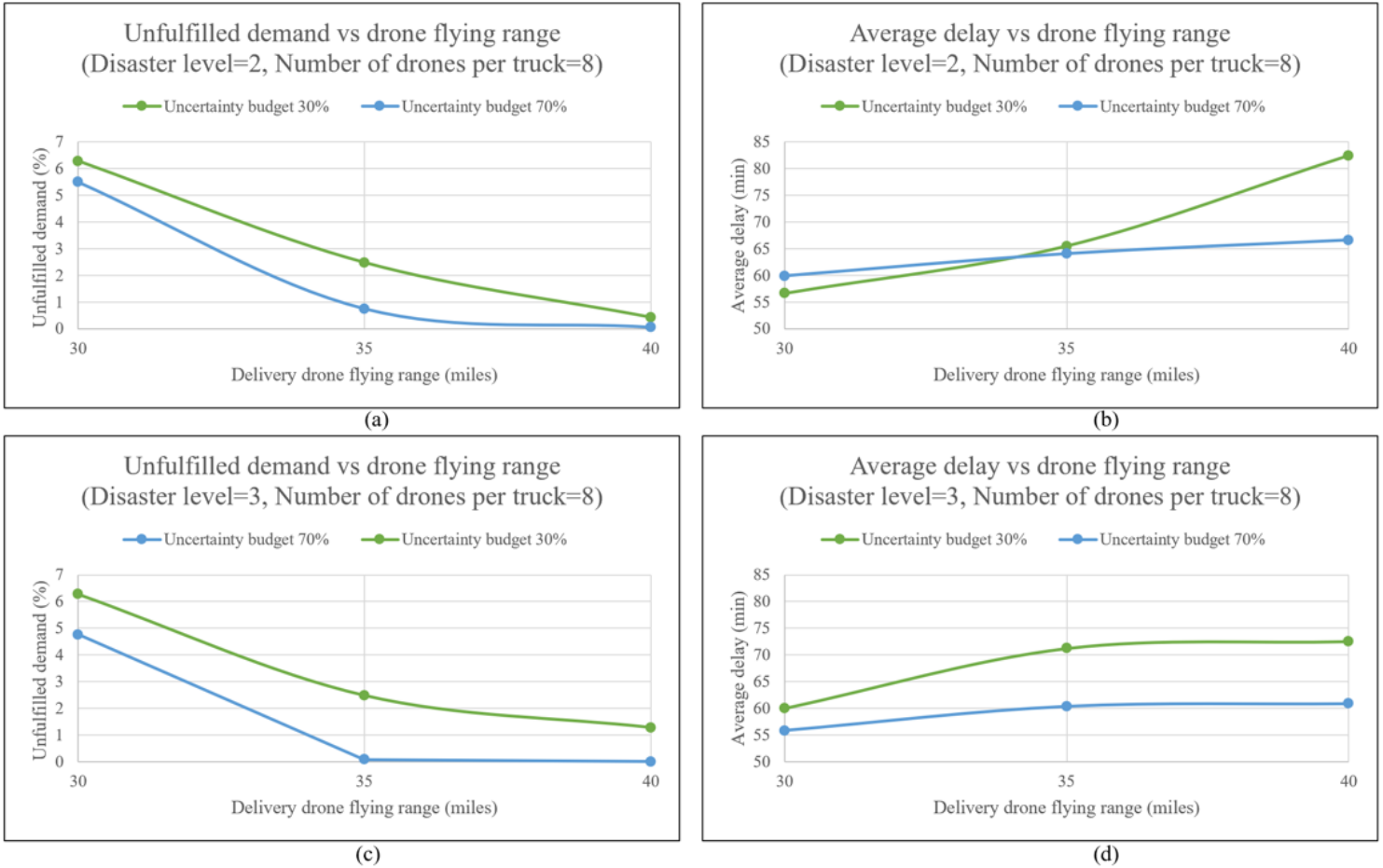}
\end{figure}

\begin{figure}[htbp]
    \centering
    \caption{Variation in average unfulfilled demand and average delay with delivery drone flying range (6 delivery drones per truck). \label{fig:range_of_drones_6}}
    \includegraphics[width=.99\textwidth]{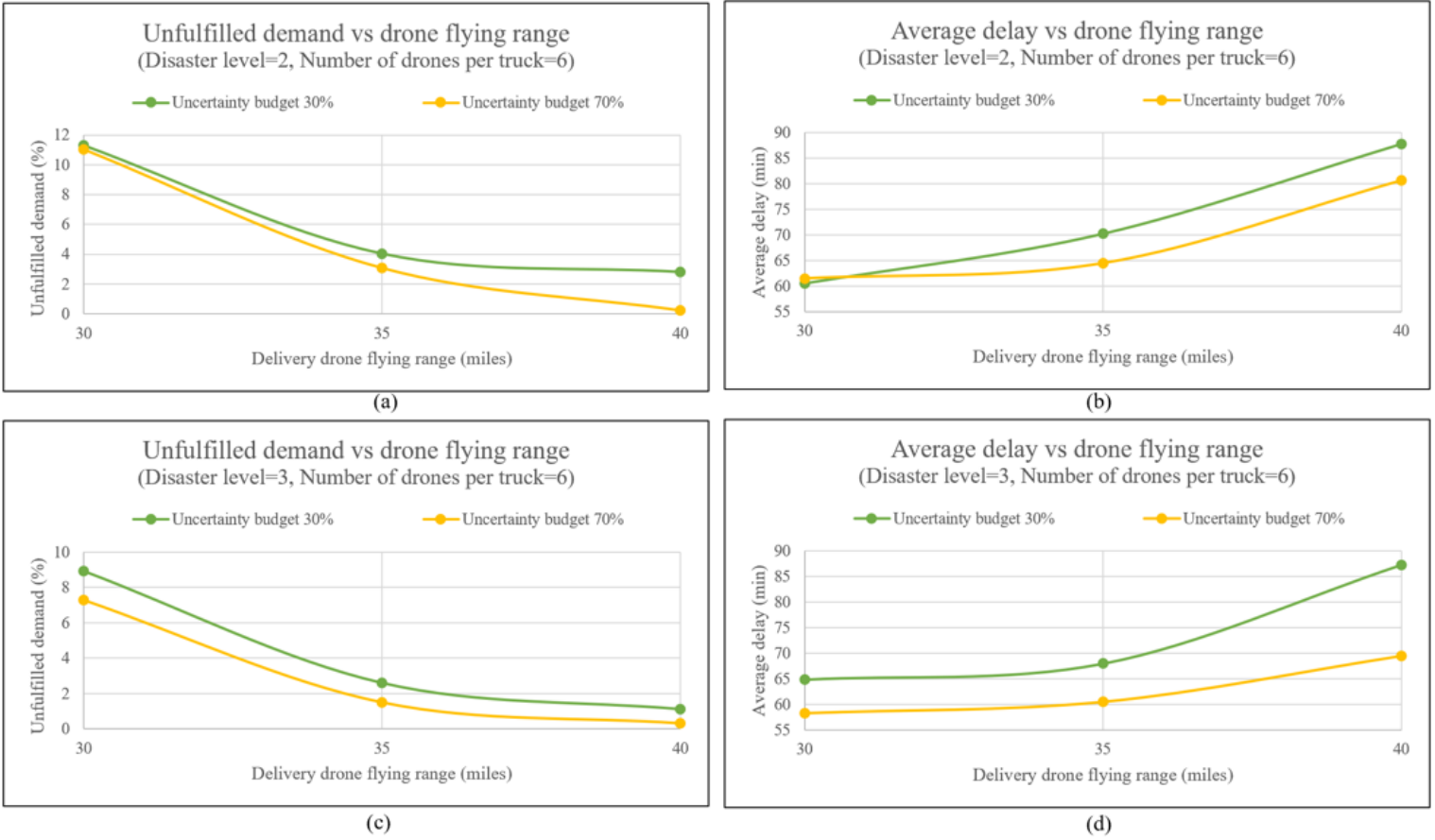}
\end{figure}

We also study the effect of drone flying range on the average proportion of unfulfilled demand and average delay time and we present our findings in Figures \ref{fig:range_of_drones_8} and \ref{fig:range_of_drones_6}. We note that when the number of delivery drones carried by each truck is fixed, the proportion of unfulfilled demand decreases as we increase the drone flying range.

When using 8 drones per truck, both disaster levels 2 and 3 (as seen in Figures \ref{fig:range_of_drones_8}(a) and \ref{fig:range_of_drones_8}(c)) reveal that the proportion of unfulfilled demand remains higher for the 30\% budget uncertainty cases as compared to the 70\% uncertainty budget cases. This is consistent with our explanations and findings in Figure \ref{fig:number_of_drones}. Comparing Figures \ref{fig:range_of_drones_8}(b) and \ref{fig:range_of_drones_8}(d), we also observe that the average delay in reaching communities decreases with the increase in drone flying range. The rate of increase is higher in the case of uncertainty budget of 30\% and disaster level 2. We also find that once the average unfulfilled proportion of demand reaches 0\%, the average delay time remains almost unchanged even when the drone flying range is increased. This is a similar finding to our observation when investigating the number of drones per truck in subsection \ref{subsec:number_of_drones}.

The difference between cases with uncertainty budget 30\% and 70\% in terms of the average proportion of unfulfilled demand is evident in Figure \ref{fig:range_of_drones_8}, where we used 8 delivery drones in each truck. When we change the number of drones per truck, and use 6 delivery drones instead of 8, the difference is made even smaller, as revealed pictorially in Figure \ref{fig:range_of_drones_6}.

In Figures \ref{fig:range_of_drones_6}(b) and \ref{fig:range_of_drones_6}(d), we see that the average delay time decreases as the flying range increases, whereas an uncertainty budget of 30\% leads to larger average delays to reach communities. However, we also find that the difference between average proportions of unfulfilled demand given by 30\% and 70\% uncertainty budgets is smaller than the case with 8 delivery drones in each truck. For both disaster level 2 in Figure \ref{fig:range_of_drones_6}(a), and disaster level 3 in Figure \ref{fig:range_of_drones_6}(c), as the average unfulfilled proportion of demand decreases, we observe the corresponding increase in the average delay in reaching communities. 

\section{Conclusions}
\label{sec:conclusions}

In this work, we study a variant of the two-echelon vehicle routing problem (2EVRP) with UAVs or drones with uncertain demand, a feature that is very useful in the context of humanitarian logistics in the aftermath of a disaster. The lack of stable infrastructure (transportation, power, communications) makes both reaching and identifying affected communities very difficult. Following a disaster, the demand for emergency supplies (including small medical kits, epi-pens, dry food and water) is uncertain in our setup. And, even if we were able to obtain exact demands, reaching populations is a hurdle when the transportation network is failing. Our work successfully addresses both these issues faced by humanitarian logistics providers. 

To handle this special case of humanitarian aid delivery problem, we propose a unique 2EVRP framework, where trucks (the first echelon vehicles) carrying drones travel along the still functional roads of the transportation network and make intermittent stops at designated locations (satellites). When a truck stops at one of those designated locations, delivery drones (the second echelon vehicles) are dispatched for aid package deliveries in the disaster affected communities. 

We develop a novel computational framework to handle the demand uncertainty in the resulting 2EVRP. In the proposed framework, we design a two-stage robust computational approach. We develop a column-and-constraint-generation approach for worst case demand scenario generation for a given truck and drone routing decision. Moreover, we also develop a CG-based decomposition scheme and heuristic to generate drone routes for a set of demand scenarios efficiently. We combine the column generation scheme within a column-and-constraint-generation approach to determine robust decisions regarding truck and drone routes, time to visit communities, and the quantities of aid materials delivered in these communities. 

We validate our framework with real-world size dataset that simulates the demand for emergency aid in different areas in Puerto Rico after a disaster. From our numerical experiments, we gather managerial insights into the impacts of various parameters on two of the most important measures of efficiency in post-disaster emergency aid delivery operations: the average proportion of unfulfilled demand and the average delay in reaching demand locations. Our framework and computational experiments will help humanitarian logistics providers in making decisions considering the geospatial impact of a disaster and the population density distribution of the disaster affected region. Humanitarian organizations can also build upon our framework to incorporate problem specific characteristics, for example, prioritizing expedient aid delivery, ensuring equity, etc. into the decision making.    

In our 2EVRP framework, we consider only one source of uncertainty, i.e., demand for emergency aid at the affected communities. During real-world humanitarian logistics operations, several other parameters are of a stochastic nature. For example, we consider the flying time range and load carrying capacities of drones as deterministic. These two properties of drones are significantly impacted by weather conditions during delivery operations and have been identified as important future work avenues in the literature (see, e.g., \cite{poikonen2021future}). We plan to extend our study by incorporating uncertain weather conditions in terms of stochastic flying range and load carrying capacity of drones in future work.



\bibliographystyle{abbrv_networks}
\bibliography{bibliography}

\newpage

\appendix

\end{document}